# STOCHASTIC CONTROL PROBLEMS FOR SYSTEMS DRIVEN BY NORMAL MARTINGALES

By Rainer Buckdahn, Jin Ma[1] and Catherine Rainer

*Université de Bretagne Occidentale, Purdue University
and Université de Bretagne Occidentale*

In this paper we study a class of stochastic control problems in which the control of the jump size is essential. Such a model is a generalized version for various applied problems ranging from optimal reinsurance selections for general insurance models to queueing theory. The main novel point of such a control problem is that by changing the jump size of the system, one essentially changes the type of the driving martingale. Such a feature does not seem to have been investigated in any existing stochastic control literature. We shall first provide a rigorous theoretical foundation for the control problem by establishing an existence result for the multidimensional structure equation on a Wiener–Poisson space, given an arbitrary bounded jump size control process; and by providing an auxiliary counterexample showing the nonuniqueness for such solutions. Based on these theoretical results, we then formulate the control problem and prove the Bellman principle, and derive the corresponding Hamilton–Jacobi–Bellman (HJB) equation, which in this case is a mixed second-order partial differential/difference equation. Finally, we prove a uniqueness result for the viscosity solution of such an HJB equation.

**1. Introduction.** In this paper we are interested in a class of stochastic control problems in which the dynamics of the controlled system take the following form:

$$(1.1) \quad Y_t = x + \int_0^t b(Y_s, u_s, \pi_s)\, ds + \int_0^t \sigma(Y_{s-}, \pi_s, u_s)\, dX_s^u, \qquad t \geq 0,$$

Received December 2006; revised June 2007.
[1]Supported in part by NSF Grant 0204332.
*AMS 2000 subject classifications.* Primary 93E20; secondary 60G44, 35K55.
*Key words and phrases.* Normal martingales, structure equation, Bellman principle, HJB equation, partial differential/difference equation, viscosity solution.

where $X^u$ is a martingale that satisfies the so-called *structure equation*:

$$(1.2) \qquad [X^u]_t = t + \int_0^t u_s \, dX^u_s, \qquad t \geq 0.$$

In the above $[X^u]$ denotes the quadratic variation of $X^u$, and $u$ is some predictable process. It is easily seen that the process $u$ "controls" exactly the jumps of $X^u$, whence that of $Y$. We note that a martingale $X^u$ satisfying the structure equation (1.2) must satisfy $\langle X^u \rangle_t = t$, that is, it is a so-called *normal martingale* (cf. Dellacherie, Maisonneuve and Meyer [8]). In fact, it is known that if a normal martingale $X^u$ has the so-called *representation property*, then it must satisfy the structure equation (see, e.g., Émery [10] or Section 2 for more detailed discussions). Typical examples of normal martingales satisfying structure equation include Brownian motion, compensated Poisson process and Azéma martingale, etc.

A stochastic control problem with control being the jump size can be seen from the following example, which more or less motivated our study. In an optimal reinsurance and investment selection problem, the dynamics of the risk reserve of the insurance company can be described, in general, by a stochastic differential equation (SDE) of the following form (see, e.g., [16]):

$$(1.3) \qquad \begin{aligned} Y_t &= y + \int_0^t b(Y_s, \alpha_s(\cdot), \pi_s) \, ds + \int_0^t \sigma(\pi_s) \, dW_s \\ &\quad - \int_0^t \int_{\mathbb{R}_+} \alpha_s(x) f(s,x) \tilde{N}(dx, ds), \end{aligned}$$

where $W$ is a standard $d$-dimensional Brownian motion, representing the uncertainty of the underlying security market; $S_t \triangleq \int_0^t \int_{\mathbb{R}_+} f(s,x) \tilde{N}(dx, ds)$, with $\tilde{N}$ being a compensated Poisson random measure, denotes the accumulated incoming claims up to time $t$; the random field $\alpha$ is the "*reinsurance policy*" (or "*retention ratio*") and $\pi = (\pi^1, \ldots, \pi^n)$ is the usual investment portfolio. If we consider only those reinsurance policies $\alpha(\cdot)$ for which there exist predictable processes $\beta$ and $u$ such that the following equation holds:

$$(1.4) \qquad X_t^{0,u} \triangleq \int_0^t \beta_s \, dW_s + \int_0^t \int_{\mathbb{R}_+} \alpha_s(x) f(s,x) \tilde{N}(dx, ds),$$

where $X^{0,u}$ satisfies the structure equation (1.2), then, noting that the Brownian motion $W$ itself satisfies the structure equation with ($u \equiv 0$), we can then rewrite equation (1.3) as the general form of (1.1), in which $X^u \triangleq (W, X^{0,u})^T$ is now a $(d+1)$-dimensional normal martingale. In fact, if the probability space is properly chosen, then one can show that (see Section 2 for details) for any bounded predictable process $u$, there are always such $\alpha$ and $\beta$.

We remark that the special nature of the normal martingale gives rise to a financial market model that is complete, but with jumps (cf. Dritschel and Protter [9]). This is based on an Ocone–Haussmann–Clark type formula for normal martingales established in Ma, Protter and San Martin [17]. We should also note that as a special class of normal martingales, the solutions to structure equations are interesting in their own right, and have been studied by many authors (see, e.g., Attal and Émery [3, 4], Taviot [21] and, very recently, Émery [11]; for the two-dimensional case, see also Attal and Émery [3, 4] and Kurtz [14]). However, this is a subject that has not been explored fully. In fact, to our best knowledge, there has not been any general result regarding the well-posedness of such equations in higher dimensional cases. In this paper we show that, in a Wiener–Poisson space, the issue becomes much more tractable. In fact, we shall derive a necessary and sufficient condition for a (multidimensional) normal martingale to satisfy a structure equation in such a space. This characterization theorem will then lead to the existence theorem of the (multidimensional) structure equation in a Wiener–Poisson space. An interesting observation, however, is that such a solution is *not* unique, even in the sense of law(!). We shall provide a counterexample in the Appendix for the interested reader.

As the first step toward the full understanding of this new type of stochastic control problem, we shall first establish the dynamic programming principle and study the corresponding HJB equation. Several technical difficulties arise immediately. For example, the nonuniqueness of the solution to the structure equation requires a careful formulation of the control problem. For this reason, we formulate the control problem over the canonical Wiener–Poisson space. Also, since the control actions actually change the type of the system (as a semimartingale), we need a general Itô formula that covers all possible cases in a unified form. Furthermore, it is necessary to validate the Bellman principle in this new setting so that all our arguments will have a rigorous theoretical foundation. It should also be noted that in this case the HJB equation takes a new form which we shall name as a "*mixed second-order partial differential/difference equation*," depending on whether or not the jump part of the control is present. To our best knowledge, such a type of HJB equation is novel. As the main results of this paper, we prove that the value function is the *unique* viscosity solution of the HJB equation.

The rest of the paper is organized as follows. In Section 2 we give necessary background theory on normal martingales and, in particular, the properties of martingales satisfying structure equations in general Wiener–Poisson spaces. In Section 3 we formulate the control problem and prove the continuity of the value function. In Section 4 we prove the Bellman principle, and in Section 5 we verify that the value function is a viscosity solution to the HJB equation. In Section 6 we prove the uniqueness of the viscosity solution. A counterexample that shows the nonuniqueness of the solution of the structure equation is given in the Appendix.

**2. Preliminaries.** Throughout this paper, we denote by $(\Omega, \mathcal{F}, P; \mathbf{F})$ a filtered probability space, and we always assume that the filtration $\mathbf{F} \stackrel{\triangle}{=} \{\mathcal{F}_t\}_{t\geq 0}$ satisfies the "*usual hypotheses*" (see, e.g., Protter [20]). The following notation will be used frequently in the sequel. Let $(U, \mathcal{B}(U), m)$ be a generic measure space, with the $\sigma$-finite measure $m$, and $\mathbb{E}$ be a general Euclidean space with Lebesgue measure. For any $1 \leq p < \infty$, we denote $L^p_{\mathbf{F}}([0,T] \times U, dt \times dm; \mathbb{E})$ to be the space of all random fields $\varphi: [0,T] \times \Omega \times U \mapsto \mathbb{E}$ such that:

(i) for fixed $u \in U$, $\varphi(\cdot, \cdot, u)$ is **F**-progressively measurable;
(ii) it holds that

$$(2.1) \qquad E \int_0^T \|\varphi(t, \cdot, \cdot)\|^p_{L^p(U,m)} \, dt < \infty.$$

When $p = 2$ and the random fields are actually **F**-progressively measurable processes, the space is denoted by $L^2_{\mathbf{F}}([0,T]; \mathbb{E})$ as usual. We shall denote $\mathcal{M}^2_0(\mathbf{F}, \mathbb{R}^d)$ to be the space of all $\mathbb{R}^d$-valued, square integrable martingales $X$ defined on $(\Omega, \mathcal{F}, P; \mathbf{F})$ such that $X_0 = 0$.

2.1. *Normal martingales and structure equations.* We first recall from [8] that a square integrable martingale $X$ is called "*normal*" if $\langle X \rangle_t = t$. Here $\langle X \rangle$ is the conditional quadratic variation process of $X$, or the compensator of the bracket process $[X]$. Since the processes $[X]$ and $\langle X \rangle$ differ by a martingale, if $X$ also has the "representation property," then it is readily seen that there exists an **F**-predictable process $u$ such that

$$(2.2) \qquad [X]_t - t = \int_0^t u_s \, dX_s \qquad \forall t \geq 0.$$

Equation (2.2) is called the *structure equation driven by $u$ and solved by $X$*, and was first studied by Émery [10]. Examples of normal martingales satisfying the structure equation include the following: Brownian motion ($u \equiv 0$), compensated Poisson process ($u \equiv \alpha \in \mathbb{R}^* \stackrel{\triangle}{=} \mathbb{R} \setminus \{0\}$), in which case $X = \alpha(N_{t/\alpha^2} - t/\alpha^2)$, where $N$ is a standard Poisson process, as well as the Azéma martingale ($u_t = -X_{t-}$) and the "parabolic" martingale ($u_t = -2X_{t-}$), etc.

The general results of existence of the solution to the structure equation have been studied by several authors, but mostly restricted to the one-dimensional case. For example, in Émery [10] it is proved that for any continuous function $f: \mathbb{R} \mapsto \mathbb{R}$ there exists at least one $X \in \mathcal{M}^2_0(\mathbf{F}; \mathbb{R})$ defined on a suitable filtered probability space $(\Omega, \mathcal{F}, P, \mathbf{F})$, which solves the structure equation with $u_t = f(X_{t-})$, $t \geq 0$. Moreover, in Émery [11] it is also shown that if $f: \mathbb{R} \mapsto \mathbb{R}$ is a more general deterministic function, and $u_t = f(t)$, $t \geq 0$, then the solution to the structure equation (2.2) exists, and it is unique in law. Other references regarding the well-posedness of structure equations can be found, for example, in Meyer [18]; Kurtz and Protter [15], Azéma and Rainer [5] or, more recently, Phan [19]. We should note that to date neither the uniqueness nor the multidimensional existence has been thoroughly explored. In fact, in general, even in a very special Wiener–Poisson space, for a given predictable process $u$, the solution to the structure equation (2.2) may not be unique, not even in law(!). We shall provide a counterexample in the Appendix to clarify this issue.

To end this subsection, we list some properties of the solution of a structure equation, which will be useful in the sequel. If $X \in \mathcal{M}_0^2(\mathbf{F}, \mathbb{R})$ is a solution to the structure equation (2.2) driven by the process $u$, we denote $\Delta X_t \stackrel{\triangle}{=} X_t - X_{t-}$, and for each $\omega \in \Omega$, let $D_X(\omega) \stackrel{\triangle}{=} \{t > 0; \Delta X_t(\omega) \neq 0\}$. Then:

(i) For $P$-a.e. $\omega \in \Omega$, it holds that $\Delta X_t = u_t$, for all $t \in D_X(\omega)$;

(ii) The continuous and the pure jump part of the martingale $X$, denoted by $X^c$ and $X^d$, satisfy respectively,

$$(2.3) \qquad dX_t^c = \mathbf{1}_{\{u_t=0\}} \, dX_t \quad \text{and} \quad dX_t^d = \mathbf{1}_{\{u_t \neq 0\}} \, dX_t, \qquad t \geq 0.$$

Finally, we recall that for any $X \in \mathcal{M}^2(\mathbf{F}, P)$, it holds that

$$(2.4) \qquad [X]_t = \sum_{0 < s \leq t} (\Delta X_s)^2 + \langle X^c \rangle_t = \lim_{|\pi| \to 0} \sum_{i=0}^{n-1} (X_{t_{i+1}} - X_{t_i})^2,$$

where $\pi : 0 = t_0 < \cdots < t_n = t$ is any partition of the interval $[0, t]$, and $|\pi| \stackrel{\triangle}{=} \sup_i |t_{i+1} - t_i|$ denotes the mesh size of $\pi$. Further, the limit in (2.4) is in probability.

2.2. *The Wiener–Poisson space.* Since the well-posedness of the structure equation is essential in a control problem, we shall first take a closer look at this issue, in a special probability space: the Wiener–Poisson space. To be more precise, let $(\Omega, \mathcal{F}, P)$ be some probability space on which is defined a $d$-dimensional standard Brownian motion $B = \{B_t : t \geq 0\}$ and a (time-homogeneous) Poisson random measure $\mu$, defined on $[0, T] \times \mathbb{R}$. We assume that $B$ and $\mu$ are independent, and that the Lévy measure of $\mu$, denoted by $\nu$, satisfies the standard integrability condition

$$(2.5) \qquad \int_{\mathbb{R}^*} (1 \wedge |x|^2) \nu(dx) < +\infty,$$

where $\mathbb{R}^* \stackrel{\triangle}{=} \mathbb{R} \setminus \{0\}$. For simplicity, we assume that $\nu(\{0\}) = 0$. We denote $\mathbf{F}^{B,\mu} = \{\mathcal{F}_t^{B,\mu}\}_{t \geq 0}$ to be the natural filtration generated by $B$ and $\mu$, that is, $\mathcal{F}_t^{B,\mu} = \sigma\{B_s, \mu([0, s] \times A); 0 \leq s \leq t, A \in \mathcal{B}(\mathbb{R}^*)\}$, $t \geq 0$, and denote by $\mathbf{F}$ the augmented version of $\mathbf{F}^{B,\mu}$. Then $\mathbf{F}$ satisfies the usual hypotheses. In the rest of the paper we shall content ourselves to this probability space without further specification.

One of the most important features for the Wiener–Poisson space defined above is the following *martingale representation theorem* (see, e.g., Jacod and Shiryaev [13]): for any $X \in \mathcal{M}^2(\mathbf{F}; \mathbb{R}^d)$ such that $X_0 = 0$, there exists a unique pair of processes $(\alpha, \beta) \in L^2_{\mathbf{F}}([0,T]; \mathbb{R}^{d \times d}) \times L^2_{\mathbf{F}}([0,T] \times \mathbb{R}^*; dt \times d\nu; \mathbb{R}^d)$ for each $T > 0$, such that

$$(2.6) \qquad X_t = \int_0^t \alpha_s \, dB_s + \int_0^t \int_{\mathbb{R}^*} \beta_s(x) \tilde{\mu}(dx, ds), \qquad t \geq 0.$$

Here $\tilde{\mu}(dt\, dx) \stackrel{\triangle}{=} \mu(dt\, dx) - \nu(dx)\, dt$, $(t,x) \in [0,\infty) \times \mathbb{R}^*$, is the compensated Poisson random measure. In what follows we call the pair $(\alpha, \beta)$ in (2.6) the "representation kernel" of $X$.

We now give a result that describes the necessary and sufficient conditions for a square-integrable martingale to satisfy the structure equation in a Wiener–Poisson space, which will play a crucial role in our future discussions.

PROPOSITION 2.1. *Let $u = \{u_t\}_{t \geq 0}$ be a bounded, $\mathbf{F}$-predictable process taking values in $\mathbb{R}^d$; and let $X \in \mathcal{M}^2_0(\mathbf{F}; \mathbb{R}^d)$ that has a representation kernel $(\alpha, \beta)$ in the sense of (2.6). Then, $X$ satisfies a structure equation*

$$(2.7) \qquad \begin{cases} d[X^i]_t = dt + u^i_t \, dX^i_t, & 1 \leq i \leq d, \\ d[X^i, X^j]_t = 0, & 1 \leq i < j \leq d, t \geq 0, \end{cases}$$

*if and only if there are random sets $A^i_s \in \mathcal{B}(\mathbb{R}^*) \otimes \mathcal{F}_s$, $s \geq 0$, $1 \leq i \leq d$, such that:*

(i) $\sum_{k=1}^d \alpha^{i,k}_s \alpha^{j,k}_s = \delta_{i,j} \mathbf{1}_{\{u^i_t = 0\}}$, $dt \times dP$-*a.e.*;
(ii) $\beta^i_t(x) = u^i_t \mathbf{1}_{A^i_t}(x)$, $dt \times d\nu \times dP$-*a.e.*;
(iii) $\nu(A^i_t \cap A^j_t) \mathbf{1}_{\{u^i_t \neq 0, u^j_t \neq 0\}} = \delta_{i,j} \frac{1}{(u^i_t)^2} \mathbf{1}_{\{u^j_t \neq 0\}}$, $dt \times dP$-*a.e.*, $1 \leq i, j \leq d$.

*In the above, "$\delta_{ij}$" is Kroneker's delta.*

PROOF. Let $X \in \mathcal{M}^2_0(\mathbf{F}; \mathbb{R}^d)$ be given and let $(\alpha, \beta) \in L^2_{\mathbf{F}}([0,T]; \mathbb{R}^{d \times d}) \otimes L^2_{\mathbf{F}}(([0,T] \times \mathbb{R}^*; dt \times d\nu); \mathbb{R}^d)$ be defined by (2.6). Then, the quadratic covariation process of $X^i$ and $X^j$, for $1 \leq i, j \leq d$, is given by

$$[X^i, X^j]_t = \sum_{k=1}^d \int_0^t \alpha^{i,k}_s \alpha^{j,k}_s \, ds + \int_0^t \int_{\mathbb{R}^*} \beta^i_s(x) \beta^j_s(x) \mu(dx, ds)$$

$$(2.8) \qquad = \int_0^t \left\{ \sum_{k=1}^d \alpha^{i,k}_s \alpha^{j,k}_s + \int_{\mathbb{R}^*} \beta^i_s(x) \beta^j_s(x) \nu(dx) \right\} ds$$

$$+ \int_0^t \int_{\mathbb{R}^*} \beta_s^i(x)\beta_s^j(x)\tilde{\mu}(dx,ds).$$

Here we should note that

(2.9)
$$E\left[\int_0^t \int_{\mathbb{R}^*} |\beta_s^i(x)\beta_s^j(x)|\mu(dx,ds)\right]$$
$$= E\left[\int_0^t \int_{\mathbb{R}^*} |\beta_s^i(x)\beta_s^j(x)|\nu(dx)\,ds\right]$$
$$\leq \left\{E\left[\int_0^t \int_{\mathbb{R}^*} |\beta_s^i(x)|^2 \nu(dx)\,ds\right]\right\}^{1/2}$$
$$\times \left\{E\left[\int_0^t \int_{\mathbb{R}^*} |\beta_s^j(x)|^2 \nu(dx)\,ds\right]\right\}^{1/2} < \infty.$$

Therefore, the last integral process in (2.8) is well defined, and it is a martingale obtained by compensating the stochastic integral $\int_0^t \int_{\mathbb{R}^*} \beta_s^i(x)\beta_s^j(x)\mu(dx,ds)$, $t \geq 0$.

Now suppose that $X$ satisfies (2.7). Then, plugging $X$ with the form (2.6) into (2.7) and comparing it with (2.8), the uniqueness of the predictable semimartingale decomposition then leads us to the following identities:

(a) $u_t^i(\omega)\alpha_t^{i,j}(\omega) = 0$, for $dt \times dP$-a.e. $(t,\omega)$; $1 \leq i,j \leq d$;
(b) $\sum_{k=1}^d \alpha_t^{i,k}(\omega)\alpha_t^{j,k}(\omega) + \int_{\mathbb{R}^*} \beta_t^i(x,\omega))\beta_t^j(x,\omega))\nu(dx) = \delta_{i,j}$, for $dt \times dP$-a.e. $(t,x,\omega)$, $1 \leq i,j \leq d$;
(c) $\beta_t^i(x,\omega)\beta_t^j(x,\omega) = \delta_{i,j} \cdot u_t^i \beta_t^i(x,\omega)$, for $dt \times d\nu \times dP$-a.e. $(t,x,\omega)$, $1 \leq i,j \leq d$.

Clearly, (a) implies that $\alpha_t^{i,j}(\omega)\mathbf{1}_{\{u_t^i(\omega) \neq 0\}} = 0$, for $dt \times dP$-a.e. $(t,\omega)$, $1 \leq i, j \leq d$. Also, if we define $A_t^i \triangleq \{(x,\omega) \in \mathbb{R}^* \times \Omega : \beta_t^i(x,\omega) \neq 0\}$, then it follows from (c) that

(2.10) $\quad \beta_t^i(x,\omega) = u_t^i(\omega)\mathbf{1}_{A_t^i}(x,\omega), \qquad dt \times d\nu \times dP$-a.e. $(t,x,\omega),\ 1 \leq i \leq d$,

proving (ii). Now, setting $i = j$ in (b) and using (2.10), we have, $dt \times dP$-almost surely,

(2.11)
$$1 = \sum_{k=1}^d |\alpha_t^{i,k}|^2 \mathbf{1}_{\{u_t^i = 0\}} + \int_{\mathbb{R}^*} |u_t^i|^2 \mathbf{1}_{A_t^i}(x)\nu(dx)$$
$$= \sum_{k=1}^d |\alpha_t^{i,k}|^2 \mathbf{1}_{\{u_t^i = 0\}} + |u_t^i|^2 \nu(A_t^i).$$

Multiplying $\mathbf{1}_{\{u_t^i = 0\}}$ and $\mathbf{1}_{\{u_t^i \neq 0\}}$ on both sides above, respectively, one can easily see that $\sum_{k=1}^d |\alpha_t^{i,k}|^2 = 1$ on $\{u_t^i = 0\}$ and $\nu(A_t^i) = \frac{1}{|u_t^i|^2}$ on $\{u_t^i \neq 0\}$, $1 \leq i \leq d$, $dt \times dP$-a.e. On the other hand, setting $i \neq j$ in (b) we have, for $dt \times dP$-a.e.,

$$
\begin{aligned}
0 &= \sum_{k=1}^{d} \alpha_t^{i,k} \alpha_t^{j,k} + \int_{\mathbb{R}^*} \beta_t^i(x) \beta_t^j(x) \nu(dx) \\
&= \sum_{k=1}^{d} \alpha_t^{i,k} \alpha_t^{j,k} \mathbf{1}_{\{u_t^i = 0, u_t^j = 0\}} + u_t^i u_t^j \nu(A_t^i \cap A_t^j).
\end{aligned}
\tag{2.12}
$$

We can then easily check that $\nu(A_t^i \cap A_t^j) = 0$ on the set $\{u_t^i \neq 0, u_t^j \neq 0\}$, hence, (iii) holds. This then further implies that $\sum_{k=1}^{d} \alpha_t^{i,k} \alpha_t^{j,k} = 0$, $dt \times dP$-a.e., proving (i).

Conversely, if $X \in \mathcal{M}_0^2(\mathbf{F}; \mathbb{R}^d)$ has the form

$$X_t = X_0 + \int_0^t \alpha_s \, dB_s + \int_0^t \int_{\mathbb{R}^*} \beta_s \tilde{\mu}(dx, ds), \qquad t \geq 0,$$

with $(\alpha, \beta)$ satisfying (i)–(iii), then a straightforward calculation starting from (2.8) shows that it satisfies the structure equation (2.7). The proof is now complete. $\square$

As a direct consequence of Proposition 2.1, we can now prove an existence result for the structure equation on the Wiener–Poisson space. We note that a Lévy measure can have at most countably many atoms, we shall focus on the continuous part of the Lévy measure. In other words, denoting the set of all atoms of $\nu$ by $\Gamma$, thus, the "continuous part" of $\nu$, denoted by $\nu^c$, is defined by

$$\nu^c(A) = \nu(A) - \sum_{z \in \Gamma \cap A} \nu(\{z\}) \qquad \forall A \in \mathcal{B}(\mathbb{R}^*), 0 \notin \overline{A},$$

where $\overline{A}$ is the closure of $A$ in $\mathbb{R}$. We have the following existence result for a (possibly) multidimensional structure equation.

THEOREM 2.2. *Assume that $\nu^c([-1, 1]) = +\infty$. Then, for any bounded, $\mathbf{F}$-predictable process $u_t = (u_t^1, \ldots, u_t^d)$, $t \geq 0$, the structure equation (2.7) has at least one solution in the Wiener–Poisson space $(\Omega, \mathcal{F}, P, \mathbf{F}, B, \mu)$.*

PROOF. For each $t \geq 0$, let us define random times $1 = \tau_t^0 > \tau_t^1 > \cdots > \tau_t^d > 0$ inductively as

$$(2.13) \quad \tau_t^i = \begin{cases} \sup\{r < \tau_t^{i-1} : \nu^c((-\tau_t^{i-1}, -r] \cup [r, \tau_t^{i-1})) = (u_t^i)^{-2}\}, & u_t^i \neq 0, \\ \tau_t^{i-1}, & u_t^i = 0. \end{cases}$$

This is always possible because $\nu^c([-1, 1]) = +\infty$. Clearly, all $\tau_t^i$'s are $\mathcal{F}_t$-measurable. Moreover, the mapping $(t, \omega) \mapsto \tau_t^i(\omega)$ is jointly measurable on $\mathbb{R}^* \times \Omega$. In light of Proposition 2.1 we see that, for a given bounded, **F**-predictable process $u = (u^1, \ldots, u^d)$, we can construct a solution $X$ to the structure equation (2.7) with the form of (2.6) by simply setting

$$A_t^i = [(-\tau_t^{i-1}, -\tau_t^i] \cup [\tau_t^i, \tau_t^{i-1})] \cap \Gamma^c, \qquad 1 \leq i \leq d,$$

(2.14)
$$\alpha_t^{i,j} = \delta_{i,j} \mathbf{1}_{\{u_t^i = 0\}}, \qquad \beta_t^i(x) = u_t^i \mathbf{1}_{A_t^i}(x),$$

$$1 \leq i, j \leq d, x \in \mathbb{R}^*, t \geq 0.$$

It is easy to check that the conditions (i)–(iii) in Proposition 2.1 are satisfied, hence, the conclusion follows. □

To conclude this section, we give an Itô formula for processes that satisfy an equation of type (2.7).

PROPOSITION 2.3. *Let $u = \{u_t; t \geq 0\}$ be a bounded $\mathbf{F}$-predictable process with values in $\mathbb{R}^d$ and $X \in \mathcal{M}_0^2(\mathbf{F}; \mathbb{R}^d)$ a solution to the associated structure equation (2.7). Then, for any function $\varphi \in C^{1,2}([0,T] \times \mathbb{R}^d)$, the following formula holds:*

(2.15)
$$\varphi(t, X_t) - \varphi(0, 0) = \sum_{i=1}^d \int_0^t \mathscr{A}_{u_s}^i[\varphi](s, X_{s-}) \, dX_s^i$$
$$+ \int_0^t \left( \frac{\partial}{\partial s} \varphi(s, X_s) + \mathscr{L}_{u_s}[\varphi](s, X_s) \right) ds,$$

*for all $t \in [0, T]$, with*

$$\mathscr{A}_u^i[\varphi](s, x) = \mathbf{1}_{\{u^i = 0\}} \frac{\partial}{\partial x_i} \varphi(s, x) + \mathbf{1}_{\{u^i \neq 0\}} \frac{\varphi(s, x + u^i e_i) - \varphi(s, x)}{u^i},$$

(2.16) $\mathscr{L}_u[\varphi](s, x) = \sum_{i=1}^d \left( \mathbf{1}_{\{u^i = 0\}} \frac{1}{2} \frac{\partial^2}{\partial x_i^2} \varphi(s, x) \right.$
$$\left. + \mathbf{1}_{\{u^i \neq 0\}} \frac{\varphi(s, x + u^i e_i) - \varphi(s, x) - u^i \frac{\partial}{\partial x^i} \varphi(s, x)}{(u^i)^2} \right),$$

*where $\{e_1, \ldots, e_d\}$ is the canonical orthonormal basis in $\mathbb{R}^d$.*

PROOF. We first apply the general Itô formula to get, for all $t \in [0, T]$,

$$\varphi(t, X_t) - \varphi(0, X_0)$$

(2.17)
$$= \int_0^t \frac{\partial}{\partial s} \varphi(s, X_s) \, ds + \sum_{i=1}^d \int_0^t \frac{\partial}{\partial x_i} \varphi(s, X_{s-}) \, dX_s^i$$

$$+ \frac{1}{2} \sum_{i,j=1}^{d} \int_0^t \frac{\partial^2}{\partial x_i \, \partial x_j} \varphi(s, X_s) \, d\langle (X^i)^c, (X^j)^c \rangle_s$$

$$+ \sum_{0 < s \leq t} \left( \varphi(s, X_s) - \varphi(s, X_{s-}) - \sum_{i=1}^d \frac{\partial}{\partial x_i} \varphi(s, X_{s-}) \Delta X_s^i \right).$$

Since $\Delta X^i \Delta X^j = \Delta[X^i, X^j] = 0$, whenever $i \neq j$, we see that at any given time at most one coordinate of $X$ can jump. Furthermore, note that the jump size of the process $X^i$ is determined by $\Delta X_s^i = u_s^i$, whenever $\Delta X_s^i \neq 0$. Thus, by a standard calculation using (2.7) and relation (2.4) for each $X^i, 1 \leq i \leq d$, we can rewrite the last sum on the right-hand side of (2.17) as

$$\sum_{i=1}^d \sum_{0 < s \leq t} \mathbf{1}_{\{u_s^i \neq 0\}} \frac{\varphi(s, X_{s-} + u_s^i e_i) - \varphi(s, X_{s-}) - \frac{\partial}{\partial x_i} \varphi(s, X_{s-}) u_s^i}{(u_s^i)^2} (\Delta X_s^i)^2$$

$$= \sum_{i=1}^d \int_0^t \mathbf{1}_{\{u_s^i \neq 0\}} \frac{\varphi(s, X_{s-} + u_s^i e_i) - \varphi(s, X_{s-}) - \frac{\partial}{\partial x_i} \varphi(s, X_{s-}) u_s^i}{(u^i)_s^2}$$

$$\times (\mathbf{1}_{\{u_s^i \neq 0\}} \, ds + u_s^i \, dX_s^i).$$

Finally, note that $[X^i, X^j] = 0$, for all $1 \leq i < j \leq d$. We see that $\langle (X^i)^c, (X^j)^c \rangle = 0$, for $i \neq j$, and the result follows easily. □

**3. The stochastic control problem.** In this section we formulate our control problem and give some preliminary results. Due to the technical subtleties involved in the solutions of the structure equation as we indicated in Section 2, we find it more convenient to work on a *canonical Wiener–Poisson space*, which we now describe.

Let $T > 0$ be an arbitrarily given finite time horizon. For any $0 \leq s < t \leq T$, we denote $\Omega_{s,t}^1$ to be the space of all continuous functions from $[s, t]$ to $\mathbb{R}^d$ starting from 0, endowed with the sup-norm. We define $\mathcal{B}_1^0 \triangleq \mathscr{B}(\Omega_{s,t}^1)$ to be the topological $\sigma$-algebra of $\Omega_{s,t}^1$, and let $P_{s,t}^1$ be the Wiener measure on $(\Omega_{s,t}^1, \mathcal{B}_1^0)$.

Next, we denote $\Omega_{s,t}^2$ to be the set of all $\overline{\mathbb{N}}$-valued measures $q$ on $([s, t] \times \mathbb{R}^*, \mathscr{B}([s, t] \times \mathbb{R}^*))$ and denote $\mathcal{B}_2^0$ to be the smallest $\sigma$-algebra to which all mappings $q \in \Omega_{s,t}^2 \mapsto q(A) \in \mathbb{Z}^+ \cup \{\infty\}$, $A \in \mathscr{B}([s, t] \times \mathbb{R}^*)$, are measurable. Let $\mu_{s,t}(q, \cdot) \triangleq q(\cdot) \in \Omega_{s,t}^2$ to be the coordinate random measure defined on $(\Omega_{s,t}^2, \mathcal{B}_2^0)$, and denote $P_{s,t}^2$ to be the probability under which $\mu_{s,t}$ is a Poisson random measure with Lévy measure $\nu$ satisfying (2.5) and that $\nu^c([-1, 1]) = +\infty$. We note that the second condition is merely technical, which will guarantee the existence of solutions to the structure equation, thanks to Theorem 2.2.

We now define, for each $0 \leq s < t \leq T$, $\Omega_{s,t} \triangleq \Omega^1_{s,t} \times \Omega^2_{s,t}$; $P_{s,t} \triangleq P^1_{s,t} \otimes P^2_{s,t}$; and $\mathcal{B}_{s,t} \triangleq \overline{\mathcal{B}^0_1 \otimes \mathcal{B}^0_2}^{P_{s,t}}$, the completion of $\mathcal{B}^0_1 \otimes \mathcal{B}^0_2$ with respect to the probability measure $P_{s,t}$. The canonical Wiener–Poisson space is then defined as $(\Omega, \mathcal{F}, P) \triangleq (\Omega_{0,T}, \mathcal{B}_{0,T}, P_{0,T})$. We shall denote the generic element of $(\Omega, \mathcal{F}, P)$ by $\omega \triangleq (\omega_1, \omega_2)$, where $\omega_i \in \Omega^i_{0,T}$, $i = 1, 2$; and we define

$$B_t(\omega) = \omega_1(t), \qquad \mu(\omega, A) = \omega_2(A)$$

$$\forall t \in [0, T], \omega \in \Omega, \ A \in \mathscr{B}([0, T] \times \mathbb{R}^*).$$

Clearly, under the probability $P$, the process $B$ is a standard Brownian motion, and $\mu$ is a Poisson random measure with Lévy measure $\nu$, and they are independent.

Finally, for $t \in [0, T]$, we shall denote $\mathbf{F}^t = (\mathcal{F}^t_s)_{s \in [0,T]}$ to be the filtration on $[t, T]$ in the following sense:

$$\mathcal{F}^t_s \triangleq \sigma\{B_r - B_t, \mu(A) : A \in \mathscr{B}([t, r] \times \mathbb{R}^*), t \leq r \leq s\}, \qquad t \leq s \leq T,$$

and $\mathcal{F}^t_s$ is a trivial $\sigma$-field for all $s \leq t$. We assume that all $\mathbf{F}^t$'s are augmented by $P$-null sets. In particular, we set $\mathbf{F} \triangleq \mathbf{F}^0$, and $\mathcal{F}_s = \mathcal{F}^0_s$, for all $s \in [0, T]$. Furthermore, we denote by $\mathcal{M}^2(\mathbf{F}; \mathbb{R}^d)$ the set of all $\mathbb{R}^d$-valued, square-integrable $\mathbf{F}$-martingales $X$ on $(\Omega, \mathcal{F}, P; \mathbf{F})$.

We are now ready to describe our control problem. Let $U_1$ and $U^d$ be two nonempty compact sets in $\mathbb{R}$ and $\mathbb{R}^d$, respectively, and set $\overline{U} = U_1 \times U^d$.

DEFINITION 3.1. For all $t \in [0, T]$, we say that a triple $(\pi, u, X)$ is a "control at time $t$" if the following properties hold:

(i) $(\pi, u)$ is a pair of $\mathbf{F}^t$-predictable processes with values in $\overline{U}$;
(ii) $X \in \mathcal{M}^2(\mathbf{F}; \mathbb{R}^d)$ satisfies

$$(3.1) \quad \begin{cases} [X^i]_t = t + \int_0^t u^i_s \, dX^i_s, & 1 \leq i \leq d, t \in [0, T], \\ [X^i, X^j]_t = 0, & 1 \leq i < j \leq d, t \in [0, T], \\ X_s - X_t \perp\!\!\!\perp \mathcal{F}_t, & \forall s \geq t. \end{cases}$$

Here "$\perp\!\!\!\perp$" stands for "independent of." We denote by $\mathcal{U}(t)$ the set of all controls at time $t$.

The following proposition gives two basic properties of $\mathcal{U}(t)$.

PROPOSITION 3.2. (i) *For any $t \in [0, T]$ and any $\overline{U}$-valued $\mathbf{F}^t$-valued predictable couple of processes $(\pi, u)$, there exists $X \in \mathcal{M}^2(\mathbf{F}; \mathbb{R}^d)$ such that $(\pi, u, X) \in \mathcal{U}(t)$;*
(ii) *For any $0 \leq t \leq t' \leq T$, it holds that $\mathcal{U}(t') \subset \mathcal{U}(t)$.*

PROOF. The proof of (i) is quite straightforward. For any $t \in [0,T]$, we first apply Theorem 2.2 to obtain a solution to the structure equation (2.7) on $[0,t]$, denote it by $^tX$. We then apply Theorem 2.2 again to find a square integrable, $\mathbf{F}^t$-martingale $\{X_s^t, s \in [t,T]\}$ that solves the structure equation on the interval $[t,T]$. We now define a process

$$X_s = \begin{cases} {}^tX_s, & s \in [0,t), \\ X_s^t - X_t^t + {}^tX_t, & s \in [t,T]. \end{cases}$$

One can easily check that $X$ satisfies (3.1).

To see (ii), first note that, for all $s \geq 0$, $\mathcal{F}_s^{t'} \subset \mathcal{F}_s^t$. Thus, if $(\pi, u)$ is $\mathbf{F}^{t'}$-predictable, it is also $\mathbf{F}^t$-predictable. Now, let $u$ be a $U^d$-valued $\mathbf{F}^{t'}$-predictable process and $X$ be the corresponding solution of (3.1). Since $\mathbf{F}^{t'}$ is trivial before $t'$, $u$ must be a.s. deterministic, and equation (3.1) restricted to $[0,t']$ has a solution $X$ of the form

$$X_s^i = \int_0^s \mathbf{1}_{\{u_r^i=0\}} \, dB_r^i + \int_0^s \int_{\mathbb{R}^*} u_r^i \mathbf{1}_{A_r^i}(x) \mu(dx, ds), \qquad s \in [0, t'],$$

where $A_r^i$ is defined as the one in the proof of Theorem 2.2. In fact, this solution is unique in law, since $u$ is deterministic on $[0,t']$. It then follows that the process $\{X_s, s \in [0,t']\}$ is of independent increments. Hence, $X_s - X_t \perp\!\!\!\perp \mathcal{F}_t$, for $s \in [t, t']$. In particular, we have $X_{t'} - X_t \perp\!\!\!\perp \mathcal{F}_t$. Finally, for $s \in (t', T]$, we write

$$X_s - X_t = (X_s - X_{t'}) + (X_{t'} - X_t).$$

Since $X_s - X_{t'}$ is independent of $\mathcal{F}_{t'}$ by definition, whence of $\mathcal{F}_t$, we conclude that $X_s - X_t$ is independent of $\mathcal{F}_t$, proving the proposition. $\square$

We now describe the main ingredients of our control problem. Given any initial data $(t,y) \in [0,T] \times \mathbb{R}^m$ and a control $a = (\pi, u, X) \in \mathcal{U}(t)$, we assume that the *controlled dynamics* $Y^{t,y,a} = Y$ satisfy the following SDE driven by the normal martingale $X$:

$$(3.2) \quad Y_s = y + \int_t^s b(Y_r, \pi_r, u_r) \, dr + \int_t^s \sigma(Y_{r-}, \pi_r, u_r) \, dX_r, \qquad s \in [t,T],$$

and we make the convention that $Y_s \equiv y$, for all $s < t$. We consider the following *cost functional* for the control problem:

$$(3.3) \qquad J(t,y;a) \triangleq E[g(Y_T^{t,y,a})], \qquad (t,y) \in [0,T] \times \mathbb{R}^m,$$

and therefore, the *value function* is given by

$$(3.4) \qquad V(t,y) = \inf_{a \in \mathcal{U}(t)} E[g(Y_T^{t,y,a})], \qquad (t,y) \in [0,T] \times \mathbb{R}^m.$$

Throughout this paper we shall make use of the following assumptions on the coefficients $b$, $\sigma$ and $g$.

(H1) The functions $b \colon \mathbb{R}^m \times \overline{U} \to \mathbb{R}^m$ and $\sigma \colon \mathbb{R}^m \times \overline{U} \to \mathbb{R}^{m \times d}$ are uniformly continuous in $(y, \pi, u)$ and Lipschitz in $y$, uniformly with respect to $(\pi, u)$.

(H2) The function $g \colon \mathbb{R}^m \to \mathbb{R}$ is bounded and continuous.

It is well known that, under the assumptions (H1) and (H2), for any given control $a \in \mathcal{U}(0)$, there exists a unique $\mathbf{F}$-adapted continuous solution to (3.2), that we denote by $Y^{t,y,a} = (Y_s^{t,y,a})_{s \geq t}$, or $Y^{t,y}$ for simplicity. If $a \in \mathcal{U}(t)$, then $Y^{t,y}$ is $\mathbf{F}^t$-adapted. A simple application of Itô's formula (Proposition 2.3) yields the following result regarding the second-order operator associated with the process $Y^{t,y}$. We state it as a ready reference, but omit the proof.

PROPOSITION 3.3. *Let $a = (\pi, u, X) \in \mathcal{U}(0)$ and $Y = Y^{t,y}$ be the unique solution of equation (3.2). Then for any $\varphi \in C^{1,2}([0,T] \times \mathbb{R}^m)$, it holds that*

$$
\begin{aligned}
&\varphi(s, Y_s) - \varphi(t, y) \\
(3.5) \quad &= \sum_{i=1}^{d} \int_t^s \mathscr{A}_{\pi_r, u_r}^i[\varphi](r, Y_{r-}) \, dX_r^i \\
&\quad + \int_t^s \left( \frac{\partial}{\partial s} \varphi(r, Y_r) + \mathscr{L}_{\pi_r, u_r}[\varphi](r, Y_r) \right) dr, \qquad s \in [t, T],
\end{aligned}
$$

*where*

$$
\begin{aligned}
&\mathscr{A}_{\pi, u}^i[\varphi](r, y) \\
&= \mathbf{1}_{\{u_r^i = 0\}} \nabla_y \varphi(r, y) \sigma^i(y, \pi, u) \\
&\quad + \mathbf{1}_{\{u^i \neq 0\}} \frac{\varphi(r, y + u^i \sigma^i(y, \pi, u)) - \varphi(r, y)}{u^i}, \\
&\mathscr{L}_{\pi, u}[\varphi](r, y) \\
&= \nabla_y \varphi(r, y) b(y, \pi, u) \\
&\quad + \sum_{i=1}^{d} (\mathbf{1}_{\{u^i = 0\}} \tfrac{1}{2} (D_{yy}^2 \varphi(s, y) \sigma^i(y, \pi, u), \sigma^i(y, \pi, u)) \\
&\qquad + \mathbf{1}_{\{u^i \neq 0\}} (\varphi(r, y + u^i \sigma^i(y, \pi, u)) - \varphi(s, y) \\
&\qquad\qquad - u^i \nabla_y \varphi(r, y) \sigma^i(y, \pi, u))/(u^i)^2).
\end{aligned}
$$

*Above, $\sigma^i$ denotes the $i$th column of the matrix $\sigma$.* □

We now give a useful result of the value function.

PROPOSITION 3.4. *Assume* (H1) *and* (H2). *Then, the mapping* $(t,y) \mapsto V(t,y)$ *is continuous on* $[0,T] \times \mathbb{R}^m$.

PROOF. It suffices to show that, for any $M \geq 1$, $V(\cdot,\cdot)$ is continuous on $[0,T] \times \overline{B_M(0)}$, where $\overline{B_M(0)}$ denotes the closed ball in $\mathbb{R}^m$ centered at 0, with radius $M$.

To this end, let $0 \leq t \leq t' \leq T$ and $y, y' \in \overline{B_M(0)}$. Recall from Proposition 3.2(ii) that $\mathcal{U}(t') \subset \mathcal{U}(t)$. One can deduce by standard arguments that, for some $C > 0$, depending only on the coefficients $b$, $\sigma$ and $g$, the following inequalities hold for all $0 \leq t \leq t' \leq T$ and $y, y' \in \mathbb{R}^m$:

$$\sup_{a \in \mathcal{U}(t)} E \left\{ \sup_{s \in [t,T]} |Y_s^{t,y}|^2 \Big| \mathcal{F}_t \right\} \leq C(1+|y|^2),$$

(3.6)
$$\sup_{a \in \mathcal{U}(t)} E \left\{ \sup_{s \in [t',T]} |Y_s^{t,y} - Y_s^{t',y'}|^2 \Big| \mathcal{F}_t \right\}$$
$$\leq C((1+|y|^2+|y'|^2)|t-t'| + |y-y'|^2),$$

where $Y^{t,y}$ and $Y^{t',y'}$ are the solutions of (3.2) starting from $(t,y)$ and $(t',y')$, respectively.

Now, for any $\varepsilon > 0$, there is some $a = (\pi, u, X) \in \mathcal{U}(t)$ depending on $(t,y)$, such that

(3.7) $$V(t,y) \geq E[g(Y_T^{t,y})] - \varepsilon.$$

Also, applying (3.6), we can find, for the given $\varepsilon > 0$, some $N = N_\varepsilon \geq 1$ such that

(3.8) $$\sup_{a \in \mathcal{U}(t)} P\{|Y_T^{t,y}| \vee |Y_T^{t',y'}| \geq N\} \leq \frac{2C(1+M^2)}{N^2} \leq \frac{\varepsilon}{2}.$$

Since $g$ is bounded and continuous, thanks to (H2), for the given $\varepsilon > 0$ and $N = N_\varepsilon > 0$, there exists $\delta = \delta(\varepsilon, N_\varepsilon) > 0$ such that $|g(z) - g(z')| \leq \varepsilon$, whenever $z, z' \in \overline{B_N(0)}$ and $|z - z'| \leq \delta$. Consequently, we have

(3.9)
$$E\{|g(Y_T^{t,y}) - g(Y_T^{t',y'})|\}$$
$$\leq \varepsilon \|g\|_\infty + E\{|g(Y_T^{t,y}) - g(Y_T^{t',y'})|\mathbf{1}_{\{|Y_T^{t,y}| \leq N, |Y_T^{t',y'}| \leq N\}}\}$$
$$\leq \varepsilon(1 + \|g\|_\infty) + 2\|g\|_\infty P\{|Y_T^{t,y} - Y_T^{t',y'}| > \delta\}.$$

Combining (3.9) and (3.6), we obtain that

$$E\{|g(Y_T^{t,y}) - g(Y_T^{t',y'})|\} \leq \varepsilon(1+\|g\|_\infty) + \frac{2}{\delta^2}\|g\|_\infty C\{(1+2M^2)|t-t'| + |y-y'|^2\}.$$

Since $V(t',y') \leq E\{g(Y_T^{t',y'})\}$ always holds, it follows that

$$V(t,y) - V(t',y')$$
(3.10)
$$\geq E\{g(Y_T^{t,y}) - g(Y_T^{t',y'})\} - \varepsilon$$
$$\geq -\varepsilon(\|g\|_\infty + 2) - \frac{2}{\delta^2}\|g\|_\infty C(1+2M^2)\{|t-t'| + |y-y'|^2\}.$$

Taking lim inf over $(t,y),(t',y') \in [0,T] \times \overline{B_M(0)}$, we obtain that

(3.11)
$$\liminf_{0 \leq t'-t \searrow 0, |y-y'| \to 0}\{V(t,y) - V(t',y')\} \geq 0.$$

It is clear that to prove the continuity of $V$ it suffices to prove the following inequality:

(3.12)
$$\limsup_{0 \leq t'-t \searrow 0, |y-y'| \to 0,}(V(t,y) - V(t',y')) \leq 0.$$

We again follow a more or less standard procedure. Namely, we first choose $a' = (\pi', u', X') \in \mathcal{U}(t')$ so that

$$V(t',y') \geq E[g(Y_T^{t',y'})] - \varepsilon.$$

Then, following the same estimates as was done for (3.10), we can derive that

$$V(t,y) - V(t',y')$$
(3.13)
$$\leq E[g(Y_T^{t,y})] - E[g(Y_T^{t',y'})] + \varepsilon$$
$$\leq \varepsilon(\|g\|_\infty + 2) + 2\|g\|_\infty C\frac{1}{\delta^2}((1+2M^2)|t-t'| + |y-y'|^2).$$

Taking lim sup for $(t,y),(t',y') \in [0,T] \times \overline{B_M(0)}$, we derive (3.12).

Combining (3.11) and (3.12), we see that, on $[0,T] \times \overline{B_M(0)}$ the limit exists and

(3.14)
$$\lim_{0 \leq t'-t \searrow 0, |y-y'| \to 0,}(V(t,y) - V(t',y')) = 0.$$

It then follows that

(3.15)
$$\lim_{|t'-t|+|y-y'| \to 0,}|V(t,y) - V(t',y')|$$
$$= \left|\lim_{0 \leq s'-s \searrow 0, |z-z'| \to 0}(V(s,z) - V(s',z'))\right| = 0,$$

where $s = t \wedge t'$, $s' = t \vee t'$, $(z,z') = (y,y')$ if $t \leq t'$ and $(z,z') = (y',y)$ otherwise. Consequently, $V$ is continuous on $[0,T] \times \overline{B_M(0)}$, and the proposition follows. $\square$

**4. The Bellman principle.** In this section we prove the following "*dynamic programming principle*" (Bellman principle) for our control problem (3.2)–(3.4).

PROPOSITION 4.1 (Dynamic programming principle). *Assume* (H1) *and* (H2). *Then, for any* $(t,y) \in [0,T] \times \mathbb{R}^m$ *and* $0 < h \leq T - t$, *it holds that*

$$(4.1) \qquad V(t,y) = \inf_{a \in \mathcal{U}(t)} E[V(t+h, Y^{t,y,a}_{t+h})].$$

PROOF. We fix $0 \leq t \leq t + h \leq T$, $y \in \mathbb{R}^m$ and $a = (\pi, u, X) \in \mathcal{U}(t)$. We denote the corresponding solution to equation (3.2) with initial data $(t,y)$ and control $a$ by $Y^{t,y}$. For each $n \geq 1$, we set $\Gamma^n = \{k 2^{-n}; k \in \mathbb{Z}^m\}$, and for each $z = k2^{-n} \in \Gamma^n$, we consider the $m$-dimensional cube $I(z) := \prod_{i=1}^{m}[(k_i - 1)2^{-n}, k_i 2^{-n})$.

Obviously, for any $\varepsilon > 0$, $n \geq 1$, and $z \in \Gamma^n$, we can find some control $a^z = (\pi^z, u^z, X^z) \in \mathcal{U}(t+h)$ such that the associated solution $Y^z \stackrel{\triangle}{=} Y^{t+h,z,a^z}$ satisfies

$$(4.2) \qquad E[g(Y^z_T)] = E[g(Y^{t+h,z,a^z}_T)] \leq V(t+h, z) + \varepsilon.$$

Then we define the new control pair:

$$(4.3) \qquad (\widehat{\pi}_s, \widehat{u}_s) = \begin{cases} (\pi_s, u_s), & \text{if } s \in (t, t+h), \\ (\pi^z_s, u^z_s), & \text{if } Y^{t,y}_{t+h} \in I(z) \text{ and } s \in [t+h, T], \end{cases}$$

and the process

$$(4.4) \quad \widehat{X}_s = \begin{cases} X_s, & \text{if } s \in [0, t+h), \\ X_{t+h} + (X^z_s - X^z_{t+h}), & \text{if } Y^{t,y}_{t+h} \in I(z) \text{ and } s \in [t+h, T], \end{cases}$$

where $z$ runs over all $\Gamma^n$. It is not hard to check that the process $\widehat{X}$ is a solution of the structure equation (2.7) driven by $\widehat{u}$; and $(\widehat{\pi}, \widehat{u}, \widehat{X})$ is a control in $\mathcal{U}(t)$. Let us denote by $\widehat{Y}$ the corresponding solution to equation (3.2) with initial data $(t,y)$. Since the processes $\{(\pi_s, u_s, X_s), s \in [t, t+h]\}$ and $\{(\pi^z_s, u^z_s, X^z_s - X^z_{t+h}), s \in [t+h, T], z \in \Gamma^n\}$, are independent, so are the solutions $\{(\widehat{Y}_s), s \in [t, t+h]\}$ and $\{(Y^z_s), s \in (t+h, T]\}$, $z \in \Gamma^n$. Furthermore, the uniqueness of the solution to SDE (3.2) implies that $\widehat{Y}_s = Y^{t,y}_s$, for all $s \in [t, t+h]$; and for $s \in [t+h, T]$, $P$-a.s. on the set $\{Y^{t,y}_{t+h} \in I(z)\}$ it holds that

$$\widehat{Y}_s - Y^z_s = (Y^{t,y}_{t+h} - z) + \int_{t+h}^{s} (b(\widehat{Y}_r, \widehat{\pi}_r, \widehat{u}_r) - b(Y^z_r, \widehat{\pi}_r, \widehat{u}_r)) \, dr$$
$$+ \int_{t+h}^{s} (\sigma(\widehat{Y}_r, \widehat{\pi}_r, \widehat{u}_r) - \sigma(Y^z_r, \widehat{\pi}_r, \widehat{u}_r)) \, d\widehat{X}_r.$$

Thus, with some standard estimates and using the fact that $\widehat{X}$ is a normal martingale, as well as Gronwall's inequality, we obtain that there is some constant $C > 0$ (depending only on the Lipschitz constants of the coefficients $b$ and $\sigma$), such that

$$(4.5) \qquad E\left\{\sup_{s\in[t+h,T]} \sum_{z\in\Gamma^n} |\widehat{Y}_s - Y_s^z|^2 \mathbf{1}_{I(z)}(Y_{t+h}^{t,y})\right\} \leq Cm4^{-n}.$$

Consequently, from (3.6) and the continuity and boundedness of the function $g$, we can deduce that, for some $n_\varepsilon \geq 1$ and all $n \geq n_\varepsilon$,

$$(4.6) \qquad E\left\{\sum_{z\in\Gamma^n} |g(\widehat{Y}_T) - g(Y_T^z)| \mathbf{1}_{I(z)}(Y_{t+h}^{t,y})\right\} \leq \varepsilon.$$

Moreover, since the $Y_T^z$'s and $Y_{t+h}^{t,y}$ are independent, for every $n \geq N(\varepsilon)$, we have

$$\begin{aligned}
E\{g(\widehat{Y}_T)\} &\leq E\left\{\sum_{z\in\Gamma^n} g(Y_T^z)\mathbf{1}_{I(z)}(Y_{t+h}^{t,y})\right\} + \varepsilon \\
&= \sum_{z\in\Gamma^n} E\{g(Y_T^z)\} P\{Y_{t+h}^{t,y} \in I(z)\} + \varepsilon \\
&\leq \sum_{z\in\Gamma^n} (V(t+h,z) + \varepsilon) P\{Y_{t+h}^{t,y} \in I(z)\} + \varepsilon \\
&= \sum_{z\in\Gamma^n} V(t+h,z) P\{Y_{t+h}^{t,y} \in I(z)\} + 2\varepsilon \\
&= E\left\{\sum_{z\in\Gamma^n} V(t+h,z)\mathbf{1}_{I(z)}(Y_{t+h}^{t,y})\right\} + 2\varepsilon.
\end{aligned}$$

Letting $n \to +\infty$ and applying Proposition 3.4, we obtain from the bounded convergence theorem that

$$V(t,y) \leq E\{g(\widehat{Y}_T)\} \leq E\{V(t+h, Y_{t+h}^{t,y})\} + 2\varepsilon.$$

Finally, since $a = (\pi, u, X) \in \mathcal{U}(t)$ is arbitrary, we conclude that

$$V(t,x) \leq \inf_{a\in\mathcal{U}(t)} E\{V(t+h, Y_{t+h}^{t,y,a})\}.$$

To prove the converse inequality, we borrow some idea of [12]. Namely, we shall split the canonical probability space into two, and patch up the desired control.

First recall that each canonical space $\Omega_{s,t}$ is the product of the canonical Wiener space and the canonical space of the Poisson random measure, thus, we can identify the probability spaces $(\Omega, \mathbf{F}, P) = (\Omega_{0,T}, \mathbf{F}_{0,T}, P_{0,T})$ with $(\Omega_{0,t+h} \times \Omega_{t+h,T}, \mathbf{F}_{0,t+h} \otimes \mathbf{F}_{t+h,T}, P_{0,t+h} \otimes P_{t+h,T})$ by the following bijection $p: \Omega_{0,T} \mapsto \Omega_{0,t+h} \times \Omega_{t+h,T}$. For $\omega = (\omega_1, \omega_2) \in \Omega_{0,T} = \Omega_{0,T}^1 \times \Omega_{0,T}^2$, we define

$$\begin{cases} \omega_{0,t+h} \triangleq (\omega_1|_{[0,t+h]}, \omega_2|_{[0,t+h]}) \in \Omega_{0,t+h}, \\ \omega_{t+h,T} \triangleq ((\omega_1 - \omega_1(t+h))|_{[t+h,T]}, \omega_2|_{[t+h,T]}) \in \Omega_{t+h,T}, \\ p(\omega) \triangleq (\omega_{0,t+h}, \omega_{t+h,T}). \end{cases}$$

Now, for $\varepsilon > 0$, let $a = (\pi, u, X) \in \mathcal{U}(t)$ such that

$$E[g(Y_T^{t,y,a})] \leq V(t, y) + \varepsilon.$$

For all $\omega_{0,t+h} \in \Omega_{0,t+h}$, we define a process on $\Omega_{t+h,T}$ by

$$(\pi(\omega_{0,t+h}, \cdot), u(\omega_{0,t+h}, \cdot))(\omega_{t+h,T}) = (\pi, u) \circ p^{-1}(\omega_{0,t+h}, \omega_{t+h,T}).$$

For $P_{0,t+h}$-almost all $\omega_{0,t+h}$, $(\pi(\omega_{0,t+h}, \cdot), u(\omega_{0,t+h}, \cdot))$ is a version of an $\mathbf{F}_{t+h,T}$-predictable process and, following the construction of Proposition 2.1, we can prove that, if $X$ solves the structure equation (2.7) driven by $u$ on $[0, T]$, then $X(\omega_{0,t+h}, \cdot) = X \circ p^{-1}(\omega_{0,t+h}, \cdot)$ solves (2.7) driven by $u(\omega_{0,t+h}, \cdot)$ on the time interval $[t+h, T]$.

We now consider the control process $a(\omega_{0,t+h}, \cdot) \triangleq (\pi(\omega_{0,t+h}, \cdot), u(\omega_{0,t+h}, \cdot), X(\omega_{0,t+h}, \cdot))$. For each $\omega_{0,t+h} \in \Omega_{0,t+h}$, we denote $\xi(\omega_{0,t+h}) = Y_{t+h}^{t,y,a}[p^{-1}(\omega_{0,t+h}, \cdot)] \in \mathbb{R}^d$. We note that $Y_{t+h}^{t,y,a}$ is $\mathcal{F}_{t+h}$-measurable, thus, $\xi$ depends only on $\omega_{0,t+h}$. Next, we denote the solution of (3.2) on $[t+h, T] \times \Omega_{t+h,T}$ with initial condition $\xi(\omega_{0,t+h})$ by $\widetilde{Y}^{t+h,\xi,a}(\omega_{0,t+h})(\cdot) \triangleq Y^{t+h,\xi(\omega_{0,t+h}),a(\omega_{0,t+h},\cdot)}$. Then, as in [12], one shows that for $P_{0,t+h}$-almost all $\omega_{0,t+h} \in \Omega_{0,t+h}$, the paths of $Y^{t,y,a}[p^{-1}(\omega_{0,t+h}, \cdot)]$ coincide with those of $\widetilde{Y}^{t+h,\xi,a}(\omega_{0,t+h})(\cdot)$, $P_{t+h,T}$-almost surely. Moreover, for $P_{0,t+h}$-almost every $\omega_{0,t+h}$, it holds that

$$E\{g(Y_T^{t,y,a})|\mathcal{F}_{t+h}^t\}(p^{-1}(\omega_{0,t+h}, \cdot))$$
$$= E_{0,t+h}\{g(Y_T^{t+h,z,a(\omega_{0,t+h},\cdot)})\}|_{z=Y_{t+h}^{t,y,a}\circ p^{-1}(\omega_{0,t+h},\cdot)},$$

$P_{t+h,T}$-almost surely. Now, for any given $\pi_0 \in \mathcal{U}_1$, we put

$$\begin{cases} \overline{\pi}(\omega_{0,t+h}, \cdot) \triangleq \pi_0 \mathbf{1}_{[0,t+h)} + \pi(\omega_{0,t+h}, \cdot) \mathbf{1}_{[t+h,T]}, \\ \overline{u}(\omega_{0,t+h}, \cdot) \triangleq u(\omega_{0,t+h}, \cdot) \mathbf{1}_{[t+h,T]}, \\ \overline{X}^{\omega_{0,t+h}} \triangleq B \mathbf{1}_{[0,t+h)} + (X(\omega_{0,t+h}, \cdot) - X_{t+h}(\omega_{0,t+h}, \cdot) + B_t) \mathbf{1}_{[t+h,T]}. \end{cases}$$

Then, clearly $\overline{a}(\omega_{0,t+h}, \cdot) = (\overline{\pi}(\omega_{0,t+h}, \cdot), \overline{u}(\omega_{0,t+h}, \cdot), \overline{X}^{\omega_{0,t+h}}) \in \mathcal{U}(t+h)$. Consequently, by definition of $V(t+h, z)$, we have for all $z \in \mathbb{R}^m$ and $P_{0,t+h}$-a.s. every $\omega_{0,t+h}$,

$$E[g(Y_T^{t+h,z,a(\omega_{0,t+h},\cdot)})] = E[g(Y_T^{t+h,z,\overline{a}(\omega_{0,t+h},\cdot)})] \geq V(t+h, z).$$

This implies that

$$V(t,y) + \varepsilon \geq E[E[g(Y_T^{t,y,a})|\mathbf{F}_{t+h}^t]]$$

(4.7)
$$\geq E[V(t+h,z)_{z=Y_T^{t,y,a}}] = E[V(t+h, Y_T^{t,y,a})]$$

$$\geq \inf_{a \in \mathcal{U}(t)} E[V(t+h, Y_T^{t,y,a})].$$

Since this is true for all $\varepsilon > 0$, we get the second inequality. This completes the proof. $\square$

**5. The HJB equation.** In this section we apply the Bellman principle of the previous section to derive the corresponding Hamilton–Jacobi–Bellman (HJB) equation for our control problem. To be more precise, we shall prove that the value function is a viscosity solution of the following fully nonlinear partial differential-difference equation (PDDE):

(5.1)
$$\begin{cases} -\dfrac{\partial}{\partial t} V(t,y) - \inf_{(\pi,u) \in \overline{U}} \mathscr{L}_{\pi,u}[V](t,y) = 0, & (t,y) \in [0,T] \times \mathbb{R}^m, \\ V(T,y) = g(y), & y \in \mathbb{R}^m, \end{cases}$$

where the second-order operator $\mathscr{L}_{\pi,u}$ is of the form (3.6),

$$\mathscr{L}_{\pi,u}[\varphi](t,y)$$
$$= \nabla_y \varphi(t,y) b(y,\pi,u)$$

(5.2)
$$+ \sum_{i=1}^{d} \{\mathbf{1}_{\{u^i=0\}} \tfrac{1}{2}(D_{yy}^2 \varphi(t,y) \sigma^i(y,\pi,u), \sigma^i(y,\pi,u))$$

$$+ \mathbf{1}_{\{u^i \neq 0\}}(\varphi(t, y + u^i \sigma^i(y,\pi,u)) - \varphi(t,y)$$

$$- u^i \nabla_y \varphi(t,y) \sigma^i(y,\pi,u))/(u^i)^2\},$$

(5.3)
$$\varphi \in C^{1,2}([0,T] \times \mathbb{R}^m).$$

Here again, $\sigma^i$ denotes the $i$th column of the matrix $\sigma$. We shall refer to equation (5.1) as the *HJB equation* in the sequel. We should note that such a second-order PDDE has not been studied systematically in the literature, therefore, in what follows we give some detailed investigation regarding the viscosity solution to such an equation. We begin by introducing the notion of viscosity solution, following the approach by Barles, Buckdahn and Pardoux in [6] (see also Alvarez and Tourin [1] and Amodori, Karlsen and La Chioma [2]).

DEFINITION 5.1. A continuous function $V:[0,T] \times \mathbb{R}^m \to \mathbb{R}$ is called a viscosity subsolution (resp. supersolution) of the PDDE (5.1) if:

(i) $V(T, y) \leq$ (resp. $\geq$) $g(y)$, $y \in \mathbb{R}^m$; and

(ii) for any $(t, y) \in [0, T) \times \mathbb{R}^m$ and $\varphi \in C^{1,2}([0, T] \times \mathbb{R}^m)$ such that $V - \varphi$ attains a local maximum (resp. minimum) at $(t, y)$, it holds that

$$-\frac{\partial}{\partial t}\varphi(t, y) - \inf_{(\pi, u) \in \overline{U}} \mathscr{L}^{\delta}_{\pi, u}[V, \varphi](t, y) \leq \text{(resp. } \geq \text{) } 0, \tag{5.4}$$

for all sufficiently small $\delta > 0$, where

$$\mathscr{L}^{\delta}_{\pi, u}[V, \varphi](t, y)$$

$$\stackrel{\triangle}{=} \nabla_y \varphi(t, y) b(y, \pi, u)$$

$$+ \sum_{i=1}^{d} \{ \mathbf{1}_{\{u^i = 0\}} \tfrac{1}{2} (D^2_{yy}\varphi(t, y)\sigma^i(y, \pi, u), \sigma^i(y, \pi, u))$$

$$+ \mathbf{1}_{\{0 < |u^i| \leq \delta\}} (\varphi(t, y + u^i \sigma^i(y, \pi, u)) - \varphi(t, y)$$

$$- u^i \nabla_y \varphi(t, y) \sigma^i(y, \pi, u))/(u^i)^2$$

$$+ \mathbf{1}_{\{|u^i| > \delta\}} (V(t, y + u^i \sigma^i(y, \pi, u)) - V(t, y) \tag{5.5}$$

$$- u^i \nabla_y \varphi(t, y) \sigma^i(y, \pi, u))/(u^i)^2 \}, \tag{5.6}$$

$$\varphi \in C^{1,2}([0, T] \times \mathbb{R}^m).$$

A function $V$ is called a viscosity solution of (5.1) if it is both a viscosity subsolution and a supersolution of (5.1).

REMARK 5.1. We note that the last two second-order difference quotients in (5.5) are designed to take away the possible singularity at $u = 0$ when $V$ is not smooth. Such an idea was also used in [6].

In the general theory of viscosity solutions one can often replace the local maximum and/or minimum in the definition above by the global ones. We shall verify that this can be done in our case as well. The proof follows largely the idea of [6], Lemma 3.3.

LEMMA 5.1. *In Definition 5.1 one can consider only those test functions $\varphi \in C^{1,2}([0, T] \times \mathbb{R}^m)$ such that $V - \varphi$ achieves a global maximum (for a viscosity subsolution) and a global minimum (for a viscosity supersolution), respectively, at $(t, y)$. Furthermore, the operator $\mathscr{L}^{\delta}_{\pi, u}[V, \varphi](t, y)$ can be replaced by $\mathscr{L}_{\pi, u}[\varphi](t, y)$ defined by (5.2).*

PROOF. We shall prove only the supersolution case. The proof of the subsolution case is similar, but easier. We leave it to the interested reader.

Let us first assume that $V$ is a supersolution of (5.1) in the sense of Definition 5.1. Fix any $(t,y) \in [0,T) \times \mathbb{R}^m$, and assume that $\varphi \in C^{1,2}([0,T] \times \mathbb{R}^m)$ is such that $V - \varphi$ attains its global minimum at $(t,y)$ [i.e., $V(s,z) - \varphi(s,z) \geq V(t,y) - \varphi(t,y)$ for all $(s,z) \in [0,T] \times \mathbb{R}^m$]. Then, for all $(\pi, u) \in \overline{U}$, and $i = 1, \ldots, m$, it holds that

$$V(t, y + u^i \sigma^i(y, \pi, u)) - V(t,y) \geq \varphi(t, y + u^i \sigma^i(y, \pi, u)) - \varphi(t,y).$$

Plugging this into (5.5), we obtain

$$\mathscr{L}^\delta_{\pi,u}[V, \varphi](t,y) \geq \mathscr{L}_{\pi,u}[\varphi](t,y), \qquad (\pi, u) \in \overline{U}.$$

Moreover, since $V$ is a supersolution in the sense of Definition 5.1, we have

(5.7)
$$-\frac{\partial}{\partial t}\varphi(t,y) - \inf_{(\pi,u)\in\overline{U}} \mathscr{L}_{\pi,u}[\varphi](t,y)$$
$$\geq -\frac{\partial}{\partial t}\varphi(t,y) - \inf_{(\pi,u)\in\overline{U}} \mathscr{L}^\delta_{\pi,u}[V,\varphi](t,y) \geq 0,$$

for all $\delta > 0$. Namely, (5.4) holds when $\mathcal{L}^\delta_{\pi,u}[V,\varphi]$ is replaced by $\mathcal{L}_{\pi,u}[\varphi]$.

We now prove the converse. Let $(t,y) \in [0,T) \times \mathbb{R}^m$ be fixed, and let $\varphi \in C^{1,2}([0,T] \times \mathbb{R}^m)$ be such that $V(s,z) - \varphi(s,z) \geq V(t,y) - \varphi(t,y)$ for all $(s,z)$ in some $\delta_0$-neighborhood of $(t,y)$ in $[0,T] \times \mathbb{R}^m$. We define a new function

$$\psi(s,z) \triangleq \varphi(s,z) - (\varphi(t,y) - V(t,y)) - \rho|(s,z) - (t,y)|^4,$$
$$(s,z) \in [0,T] \times \mathbb{R}^m,$$

for some $\rho > 0$ sufficiently small. By changing this function outside the $\delta_0/2$-neighborhood of $(t,y)$ if necessary, we can assume without loss of generality that $\psi(t,y) = V(t,y)$, and

$$V(s,z) \geq \psi(s,z) + \rho|(s,z) - (t,y)|^4 \qquad \forall (s,z) \in [0,T] \times \mathbb{R}^m.$$

Noting that $(\frac{\partial}{\partial t}, \nabla_y, D^2_{yy})\psi(t,y) = (\frac{\partial}{\partial t}, \nabla_y, D^2_{yy})\varphi(t,y)$, we have

$$|\mathscr{L}^\delta_{\pi,u}[V,\psi](t,y) - \mathscr{L}^\delta_{\pi,u}[V,\varphi](t,y)|$$

(5.8)
$$\leq \sum_{i=1}^d \mathbf{1}_{\{0<|u^i|\leq\delta\}} \frac{1}{(u^i)^2} |(\psi(t, y+u^i\sigma^i) - \psi(t,y))$$
$$- (\varphi(t, y+u^i\sigma^i) - \varphi(t,y))|$$
$$\leq C^4_\sigma (1+|y|)^4 d\delta^2 \rho,$$

whenever $\delta \leq \frac{\delta_0}{2C_\sigma(1+|y|)}$, where $C_\sigma > 0$ is such that $|\sigma| = |\sigma(y,\pi,u)| \leq C_\sigma(1+|y|)$, $(\pi,u) \in \overline{U}$.

Next, let $\chi \in C^\infty([0,T] \times \mathbb{R}^m)$ be a nonnegative function with $\mathrm{supp}(\chi) \subset \{(s,z) \in [0,T] \times \mathbb{R}^m : |(s,z) - (T/2, 0)| < 1\}$. Since $V$ is continuous, for any $\ell \geq 1$, we can find $\nu_\ell \in (0, \frac{1}{4\ell})$ such that, for all $(s,z), (s',z') \in [0,T] \times \mathbb{R}^m$ with

$$|(s,z) - (t,y)| \vee |(s',z') - (t,y)| \leq \ell + 1, \qquad |(s,z) - (s',z')| \leq \nu_\ell,$$

it holds that

$$|(V - \psi)(s,z) - (V - \psi)(s',z')| \leq \rho \frac{1}{(2\ell)^4}.$$

We now set $\chi_\ell(s,z) = \frac{1}{\nu_\ell^{m+1}} \chi(\frac{1}{\nu_\ell}(s,z))$, $(s,z) \in [0,T] \times \mathbb{R}^m$, and define

$$\begin{aligned}
\psi_\ell(s,z) = \psi(s,z) + \int_{[0,T] \times \mathbb{R}^m} & \left( V(s',z') - \psi(s',z') - \rho \frac{1}{(2\ell)^4} \right) \\
& \times \mathbf{1}_{\{|(s',z') - (t,y)| \in [\frac{1}{2\ell}, \ell+1]\}} \\
& \times \chi_\ell((s,z) - (s',z')) \, ds' \, dz'.
\end{aligned} \tag{5.9}$$

It can be easily verified that the sequence $\{\psi_\ell\}_{\ell \geq 1} \subset C^\infty([0,T] \times \mathbb{R}^m)$ enjoys the following properties:

(i) $\psi_\ell(s,z) = \psi(s,z)$, for all $(s,z) \in [0,T] \times \mathbb{R}^m$ with $|(s,z) - (t,y)| \notin (\frac{1}{4\ell}, \ell + 2)$;

(ii) $V(s,z) > \psi_\ell(s,z) \geq \psi(s,z)$, for all $(s,z) \in [0,T] \times \mathbb{R}^m$ with $\frac{1}{\ell} < |(s,z) - (t,y)| < \ell$;

(iii) $\psi_\ell(s,z) \geq V(s,z) - \rho\ell^{-4}$, for all $(s,z) \in [0,T] \times \mathbb{R}^m$ with $|(s,z) - (t,y)| \in [\frac{1}{\ell}, \ell]$.

We observe that these properties imply, in particular, the uniform convergence of $\psi_\ell$ to $V$ on all compacts in $[0,T] \times \mathbb{R}^m$. Moreover, due to our assumption, the fact that $V - \psi_\ell$ achieves a global minimum at $(t,y)$ implies that, for all $\ell \geq 1$,

$$-\frac{\partial}{\partial t} \psi_\ell(t,y) - \inf_{(\pi,u) \in \overline{U}} \mathscr{L}_{\pi,u} \psi_\ell(t,y) \geq 0.$$

Hence, we can find for every $\ell \geq 1$ a control state $(\pi_\ell, u_\ell) \in \overline{U}$ such that

$$-\frac{\partial}{\partial t} \psi_\ell(t,y) - \mathscr{L}_{\pi_\ell, u_\ell} \psi_\ell(t,y) \geq -1/\ell. \tag{5.10}$$

By extracting a subsequence if necessary, we can assume without loss of generality that the sequence $\{(\pi_\ell, u_\ell)\}_{\ell \geq 1}$ converges to some $(\pi, u) \in \overline{U}$ (recall that $\overline{U}$ is compact). Now,

$$\mathscr{L}_{\pi_\ell, u_\ell}[\psi_\ell](t,y) = \nabla_y \psi_\ell(t,y) b(y, \pi_\ell, u_\ell) + \sum_{i=1}^d \Theta_\ell^i, \tag{5.11}$$

where

$$\Theta_\ell^i \triangleq \mathbf{1}_{\{u_\ell^i=0\}} \frac{1}{2}(D_{yy}^2 \psi(t,y)\sigma^i(y,\pi_\ell,u_\ell), \sigma^i(y,\pi_\ell,u_\ell))$$

$$+ \mathbf{1}_{\{0<|u_\ell^i|\leq\delta\}}(\psi_\ell(t,y+u_\ell^i\sigma^i(y,\pi_\ell,u_\ell)) - \psi(t,y)$$

(5.12)
$$- u_\ell^i \nabla_y \psi(t,y)\sigma^i(y,\pi_\ell,u_\ell))/(u_\ell^i)^2$$

$$+ \mathbf{1}_{\{|u_\ell^i|>\delta\}}(\psi_\ell(t,y+u_\ell^i\sigma^i(y,\pi_\ell,u_\ell)) - \psi(t,y)$$

$$- u_\ell^i \nabla_y \psi(t,y)\sigma^i(y,\pi_\ell,u_\ell))/(u_\ell^i)^2.$$

Now, since $\psi_\ell$ converges to $V$ uniformly on compacts, it follows that

$$\lim_{\ell\to\infty} \Theta_\ell^i = \begin{cases} \dfrac{V(t,y+u^i\sigma^i(y,\pi,u)) - V(t,y) - u^i\nabla_y\psi(t,y)\sigma^i(y,\pi,u)}{(u^i)^2}, \\ \quad |u^i| > 0; \\ \dfrac{1}{2}(D_{yy}^2\psi(t,y)\sigma^i(y,\pi,u), \sigma^i(y,\pi,u)), \\ \quad \text{if there is an infinite subsequence of } u_l^i \text{ with } u_l^i = 0. \end{cases}$$

Note that the term in (5.12) that involves $\mathbf{1}_{\{0<|u^i|\leq\delta\}}$ does not necessarily have a limit, but since $\psi_\ell \geq \psi$, we can at least conclude that, for all $1 \leq i \leq d$,

$$\liminf_{\ell\to+\infty} \Theta_\ell^i \geq \mathbf{1}_{\{u^i=0\}}\frac{1}{2}(D_{yy}^2\psi(t,y)\sigma^i(y,\pi,u), \sigma^i(y,\pi,u))$$

$$+ \mathbf{1}_{\{0<|u^i|\leq\delta\}}(\psi(t,y+u^i\sigma^i(y,\pi,u)) - \psi(t,y)$$

(5.13)
$$- u^i\nabla_y\psi(t,y)\sigma^i(y,\pi,u))/(u^i)^2$$

$$+ \mathbf{1}_{\{|u^i|>\delta\}}(V(t,y+u^i\sigma^i(y,\pi,u)) - V(t,y)$$

$$- u^i\nabla_y\psi(t,y)\sigma^i(y,\pi,u))/(u^i)^2.$$

Thus, in light of (5.11) we obtain that $\liminf_{\ell\to+\infty}\mathscr{L}_{\pi_\ell,u_\ell}\psi_\ell(t,y) \geq \mathscr{L}_{\pi,u}^\delta[V,\psi](t,y)$. Combining with (5.10) and (5.8), we have

$$0 \leq -\frac{\partial}{\partial t}\psi(t,y) - \mathscr{L}_{\pi,u}^\delta[V,\psi](t,y)$$

$$\leq -\frac{\partial}{\partial t}\varphi(t,y) - \mathscr{L}_{\pi,u}^\delta[V,\varphi](t,y) + C_\sigma^4(1+|y|)^4 d\rho\delta^2,$$

whenever $\delta \leq \frac{\delta_0}{2C_\sigma(1+|y|)}$. Therefore, first letting $\rho \to 0$ and then taking the infimum over all control states $(\pi,u) \in \overline{U}$, we see that $V$ is a viscosity supersolution. $\square$

Our main result of this section is the following:

THEOREM 5.2. *The value function $V(t,y)$ defined by (3.4) is a viscosity solution of (5.1).*

PROOF. We first show that $V(\cdot,\cdot)$ is a subsolution. Let $(t,y) \in [0,T] \times \mathbb{R}^m$ be fixed. Given an arbitrary deterministic constant control $(\pi, u) \in \overline{\mathcal{U}}$, we consider $a = (\pi, u, X) \in \mathcal{U}(t)$ and denote by $Y^{t,y}$ the corresponding solution to SDE (3.2), as usual. From Proposition 4.1 we see that, for all $h > 0$ with $t + h \leq T$, it holds that

$$V(t,y) \leq E^\mu\{V(t+h, Y^{t,y}_{t+h})\}. \tag{5.14}$$

Now let $\varphi \in C^{1,2}([0,T] \times \mathbb{R}^m)$ be such tat $V - \varphi$ achieves a global maximum at $(t,y)$. Then, by Itô's formula, one has

$$0 \leq E^\mu\{V(t+h, Y^{t,y}_{t+h}) - V(t,y)\} \leq E^\mu\{\varphi(t+h, Y^{t,y}_{t+h}) - \varphi(t,y)\}$$

$$= E^\mu\left\{\int_t^{t+h} (\partial_s \varphi(s, Y^{t,y}_s) + \mathscr{L}_{\pi,u}[\varphi](s, Y^{t,y}_s))\,ds\right\}$$

$$= \int_t^{t+h} E^\mu\{\partial_s \varphi(s, Y^{t,y}_s) + \mathscr{L}_{\pi,u}[\varphi](s, Y^{t,y}_s)\}\,ds.$$

Since the process $Y^{t,y}$ has right-continuous trajectories, we see that the mapping $s \mapsto E^\mu\{\partial_s \varphi(s, Y^{t,y}_s) + \mathscr{L}_{\pi,u}[\varphi](s, Y^{t,y}_s)\}$ is right-continuous. Hence, dividing both sides of the above inequality by $h$ and taking the limit as $h$ tends to zero, we obtain that

$$-\frac{\partial}{\partial t}\varphi(t,y) - \mathscr{L}_{\pi,u}\varphi(t,y) \leq 0 \qquad \forall (\pi,u) \in \overline{U}.$$

It then follows from Lemma 5.1 that $V$ is a viscosity subsolution.

We now prove that $V(\cdot, \cdot)$ is a supersolution. Again, let $(t,y) \in [0,T) \times \mathbb{R}^m$ and $\varphi \in C^{1,2}([0,T) \times \mathbb{R}^m)$ be such that $V - \varphi$ attains a global minimum at $(t,y)$. Since the value of $\mathscr{L}_{\pi,u}\varphi(t,y)$ depends only on the values of $\varphi$ in a $C(1 + |y|)$-neighborhood of $(t,y)$, which does not depend on $(\pi, u) \in \overline{U}$, we can change the values of $\varphi$ outside of this neighborhood without changing $\inf_{\pi,u} \mathscr{L}_{\pi,u}\varphi(t,y)$. Further, since $V$ is bounded, we can also assume without loss of generality that the functions $\varphi$, as well as all its derivatives are bounded and such that, for some constant $C$, $|\psi(s,z)| \leq C(1+|z|)^{-2}$, $(s,z) \in [0,T] \times \mathbb{R}^m$, with $\psi = \partial_s \varphi, \nabla_y \varphi, D^2_{yy}\varphi$, respectively. In particular, we see that, for some $C > 0$,

$$|\mathscr{L}_{\pi,u}[\varphi](s,z)| \leq C, \qquad (s,z) \in [0,T] \times \mathbb{R}^m, \qquad (\pi,u) \in \overline{U}.$$

On the other hand, from the continuity of the mapping $(s,z,\pi,u) \to \frac{\partial}{\partial s}\varphi(s,z) + \mathscr{L}_{\pi,u}[\varphi](s,z)$, it follows that for an arbitrarily given $\varepsilon > 0$ there is some $\delta > 0$ such that

$$|(\partial_s \varphi(s,z) + \mathscr{L}_{\pi,u}[\varphi](s,z)) - (\partial_s \varphi(t,y) + \mathscr{L}_{\pi,u}[\varphi](t,y))| \leq \varepsilon, \tag{5.15}$$

for all $(s,z) \in [0,T] \times \mathbb{R}^m$ with $|(s,z) - (t,y)| \leq 2\delta$, and for all $(\pi, u) \in \overline{U}$.

Let $h \in (0, \delta)$. By the Bellman principle (Proposition 4.1), for any $\varepsilon > 0$, we can find a control $a^\varepsilon = (\pi^\varepsilon, u^\varepsilon, X^\varepsilon) \in \mathcal{U}(t)$ such that the associated dynamics $Y^{\varepsilon,t,y}$ satisfy that

(5.16) $$V(t,y) + h\varepsilon \geq E[V(t+h, Y^{\varepsilon,t,y}_{t+h})].$$

For notational simplicity, in what follows we shall drop the superscript "$\varepsilon$" from each element of the control $a^\varepsilon$, when the context is clear. Thus, following the same argument as in the subsolution case, one can show that

(5.17) $$\varepsilon h \geq E\bigg\{\int_t^{t+h} (\partial_s \varphi(s, Y^{t,y}_s) + \mathscr{L}_{\pi_s, u_s}[\varphi](s, Y^{t,y}_s))\, ds\bigg\}.$$

Moreover, if we denote $A_{\varepsilon,\delta,h} \stackrel{\triangle}{=} \{\sup_{s \in [t,t+h]} |Y^{t,y}_s - y| \geq \delta\}$ for the given constants $\varepsilon, \delta, h > 0$, then we can find some constant $C > 0$ such that

(5.18) $$\begin{aligned} P\{A_{\varepsilon,\delta,h}\} &= \frac{8}{\delta^2} E\bigg\{\sum_{i=1}^d \int_t^{t+h} |\sigma^i(Y^{t,y}_s, \pi_s, u_s)|^2\, d[X^i]_s\bigg\} \\ &\quad + \frac{2h}{\delta^2} E\bigg\{\int_t^{t+h} |b(Y^{t,y}_s, \pi_s, u_s)|^2\, ds\bigg\} \\ &= \frac{8}{\delta^2} E\bigg\{\sum_{i=1}^d \int_t^{t+h} |\sigma^i(Y^{t,y}_s, \pi_s, u_s)|^2\, ds\bigg\} \\ &\quad + \frac{2h}{\delta^2} E\bigg\{\int_t^{t+h} |b(Y^{t,y}_s, \pi_s, u_s)|^2\, ds\bigg\} \\ &\leq C(1 + |y|^2)\frac{1}{\delta^2} h. \end{aligned}$$

This, together with (5.15) and (5.17), yields, for $0 < h < \delta$,

$$\begin{aligned} \varepsilon h &\geq E^\mu \bigg\{\int_t^{t+h} (\partial_s \varphi + \mathscr{L}_{\pi_s, u_s}[\varphi])(s, Y^{t,y}_s)\, ds\, \mathbf{1}_{A^c_{\varepsilon,\delta,h}}\bigg\} - ChP\{A_{\varepsilon,\delta,h}\} \\ &\geq h\bigg(\partial_s \varphi(t,y) + \inf_{(\pi,u) \in \overline{U}} \mathscr{L}_{\pi,u}[\varphi](t,y)\bigg) - h\varepsilon - C(1+|y|^2)\frac{1}{\delta^2} h^2. \end{aligned}$$

In other words, it holds that

(5.19) $$\frac{\partial}{\partial t}\varphi(t,y) + \inf_{(\pi,u) \in \overline{U}} \mathscr{L}_{\pi,u}\varphi(t,y) \leq 2\varepsilon + C(1+|y|^2)\frac{1}{\delta^2} h.$$

Finally, first letting $h \to 0$ and then $\varepsilon \to 0$ in (5.19), we derive the desired inequality, hence, the result follows. $\square$

**6. The uniqueness of the viscosity solution.** In this section we discuss the uniqueness issue regarding the Hamilton–Jacobi–Bellman equation (5.1). For notational simplicity, we assume that $d = 1$. We should note that, given the special form of the HJB equation, its Hamiltonian, defined by

$$
\begin{aligned}
\mathcal{H}(t, y, v, p, S) &\stackrel{\triangle}{=} \inf_{(\pi, u)} \mathcal{L}_{\pi, u}[\varphi](t, y) \\
&= \inf_{(\pi, u)} \Bigg( pb(y, \pi, u) \\
&\quad + \sum_{i=1}^{d} \{ \mathbf{1}_{\{u^i = 0\}} \tfrac{1}{2}(S\sigma^i(y, \pi, u), \sigma^i(y, \pi, u)) \\
&\quad + \mathbf{1}_{\{u^i \neq 0\}} (\varphi(t, y + u^i \sigma^i(y, \pi, u)) - \varphi(t, y) \\
&\quad\quad - u^i p \sigma^i(y, \pi, u))/(u^i)^2 \} \Bigg),
\end{aligned}
$$
(6.1)

is not continuous in the variables $(p, S)$. This in fact causes some fundamental difficulties in the uniqueness proof. We shall nevertheless prove a uniqueness result under the following extra assumption on the control set $U$.

(H3) There exists a compact set $U_1 \subseteq \mathbb{R}$ such that:

(i) $0 \notin U_1$;
(ii) $U = U_1$ or $U = \{0\} \cup U_1$.

REMARK 6.1. The assumption (H3) amounts to saying that there exist positive constants $0 < \delta_0 \leq C$, such that every admissible control $u$ satisfies

$$u_t \neq 0 \quad \Longrightarrow \quad 0 < \delta_0 \leq |u_t| \leq C \qquad \forall t \in [0, T], P\text{-a.s.}$$

Such a restriction is not unusual. For example, in the insurance applications the lower bound $c$ could be thought of as the "*deductible*," while the upper bound $C$ the "*benefit limit*," of an insurance policy.

Our main result of this section is the following theorem.

THEOREM 6.2. *Assume* (H1)–(H3). *Then the value function* $V : [0, T] \times \mathbb{R}^m \to \mathbb{R}$ *defined by (3.4) is the unique viscosity solution of (5.1) among all bounded, continuous functions.*

PROOF. We begin with a slight reduction. Let $\gamma > 0$ be any given constant, and for any function $V(t,x)$, define $\hat{V}(t,x) \triangleq e^{\gamma(T-t)} V(T-t,x)$. Then, it is easy to check that $V$ is a viscosity solution of (5.1) if and only if $\hat{V}$ is a viscosity solution of the following equation:

(6.2)
$$\begin{cases} \dfrac{\partial}{\partial t} \hat{V}(t,x) + \gamma \hat{V}(t,x) - \inf_{(\pi,u) \in \overline{U}} \mathcal{L}_{\pi,u} \hat{V}(t,x) = 0, \\ \qquad\qquad\qquad\qquad\qquad (t,x) \in [0,T] \times \mathbb{R}^m, \\ \hat{V}(0,x) = e^{\gamma T} g(x), \qquad x \in \mathbb{R}^m. \end{cases}$$

We now consider two functions $V, W \in C_b([0,T] \times \mathbb{R})$, with $V$ being the subsolution of (6.2) and $W$ the supersolution of (6.2). As usual, we shall prove the uniqueness of the solution of (6.2) [hence, of (5.1)] by showing that $V \leq W$. To this end, let us suppose that

(6.3) $$\theta \triangleq \sup_{(t,x) \in [0,T] \times \mathbb{R}} (V(t,x) - W(t,x)) > 0,$$

and look for a contradiction.

For any $\varepsilon, \alpha > 0$, we consider the auxiliary function

$$\Psi_{\varepsilon,\alpha}(t,x,s,y) = V(t,x) - W(s,y) - \frac{\alpha}{2}(|x|^2 + |y|^2)$$
$$- \frac{\alpha}{2}\left(\frac{1}{T-t} + \frac{1}{T-s}\right) - \frac{1}{2\varepsilon}|x-y|^2 - \frac{1}{2\varepsilon}|s-t|^2.$$

Since $\Psi$ tends to $-\infty$ as $t,s \to T$ or $x,y \to \infty$, one concludes that for any $\varepsilon$ and $\alpha$ there exists $(\hat{t}, \hat{x}, \hat{s}, \hat{y}) = (\hat{t}_{\varepsilon,\alpha}, \hat{x}_{\varepsilon,\alpha}, \hat{s}_{\varepsilon,\alpha}, \hat{y}_{\varepsilon,\alpha}) \in ([0,T) \times \mathbb{R})^2$, such that

(6.4) $$\Psi_{\varepsilon,\alpha}(\hat{t},\hat{x},\hat{s},\hat{y}) = \max_{([0,T] \times \mathbb{R})^2} \Psi_{\varepsilon,\alpha}(t,x,s,y) \triangleq M_{\varepsilon,\alpha}.$$

Further, by (6.3) and the definition of the supremum, we see that, for all $\eta > 0$, there exists a pair $(t_\eta, x_\eta)$, $0 < t_\eta < T, x_\eta \in \mathbb{R}$, such that

(6.5) $$V(t_\eta, x_\eta) - W(t_\eta, x_\eta) \geq \theta - \eta/2.$$

Combining (6.4) and (6.5), we can find $\alpha_\eta$ such that, for all $\alpha \in (0, \alpha_\eta)$, it holds that

(6.6)
$$M_{\varepsilon,\alpha} \geq V(t_\eta, x_\eta) - W(t_\eta, x_\eta) - \alpha |x_\eta|^2 - \alpha \frac{1}{T-t_\eta}$$
$$\geq (\theta - \eta/2) - \alpha\left(|x_\eta|^2 + \frac{1}{T-t_\eta}\right) \geq \theta - \eta.$$

In other words, if $\alpha < \alpha_\eta$, then one must have

(6.7)
$$V(\hat{t}, \hat{x}) - W(\hat{s}, \hat{y}) - \frac{|\hat{x}-\hat{y}|^2 + |\hat{s}-\hat{t}|^2}{2\varepsilon}$$
$$- \frac{\alpha}{2}(|\hat{x}|^2 + |\hat{y}|^2) - \frac{\alpha}{2}\left(\frac{1}{T-\hat{t}} - \frac{1}{T-\hat{s}}\right) \geq \theta - \eta.$$

Now, for fixed $\alpha > 0$, letting $\varepsilon \to 0$, we see that

(6.8) $\qquad |\hat{x} - \hat{y}|^2 + |\hat{s} - \hat{t}|^2 \to 0, \qquad |\hat{x}|^2 + |\hat{y}|^2 \leq C_\alpha.$

Therefore, possibly along a subsequence, still denoted by $(\hat{t}, \hat{x}, \hat{s}, \hat{y})$, the family $(\hat{t}, \hat{x}, \hat{s}, \hat{y})$ converges to $(t_\alpha, x_\alpha, t_\alpha, x_\alpha)$ for some $(t_\alpha, x_\alpha) \in [0, T] \times \mathbb{R}$. Consequently, we have

$$V(\hat{t}, \hat{x}) - W(\hat{s}, \hat{y}) \to V(t_\alpha, x_\alpha) - W(t_\alpha, x_\alpha) \geq \theta.$$

Here, the last inequality is due to the definition of $\theta$ [cf. (6.3)] and that

(6.9) $\qquad \limsup_{\varepsilon \to 0} \left( \frac{1}{2\varepsilon}(|\hat{x} - \hat{y}|^2 + |\hat{s} - \hat{t}|^2) + \frac{\alpha}{2}(|\hat{x}|^2 + |\hat{y}|^2) \right) \leq \eta.$

We can now apply Ishii's lemma [7] to obtain that, for all $\rho > 0$, there exists $(X, Y) \in \mathbb{R}^{2m}$ such that

(6.10) $\qquad \begin{cases} \left( \dfrac{\hat{t} - \hat{s}}{\varepsilon} + \dfrac{\alpha}{2} \dfrac{1}{(T - \hat{t})^2}, \dfrac{\hat{x} - \hat{y}}{\varepsilon} + \alpha \hat{x}, X \right) \in \overline{\mathcal{P}}^{1,2,+} V(\hat{t}, \hat{x}), \\ \left( \dfrac{\hat{t} - \hat{s}}{\varepsilon} - \dfrac{\alpha}{2} \dfrac{1}{(T - \hat{s})^2}, \dfrac{\hat{x} - \hat{y}}{\varepsilon} - \alpha \hat{y}, Y \right) \in \overline{\mathcal{P}}^{1,2,-} W(\hat{s}, \hat{y}), \end{cases}$

and

$$\begin{pmatrix} X & 0 \\ 0 & -Y \end{pmatrix} \leq A + \rho A^2 \qquad \text{with } A = \frac{1}{\varepsilon} \begin{pmatrix} I_m & -I_m \\ -I_m & I_m \end{pmatrix} + \alpha \begin{pmatrix} I_m & 0 \\ 0 & I_m \end{pmatrix},$$

where $\overline{\mathcal{P}}^{1,2,+} V(\hat{t}, \hat{x})$ [resp. $\overline{\mathcal{P}}^{1,2,-} W(\hat{t}, \hat{x})$] denotes the "parabolic superjet" (resp. "subjets"), as defined in [7]. Choosing $\rho = \min(\varepsilon, \frac{1}{\alpha}) \leq \frac{2}{\varepsilon}$ yields

(6.11) $\qquad A + \rho A^2 \leq \dfrac{3}{\varepsilon} \begin{pmatrix} I_m & -I_m \\ -I_m & I_m \end{pmatrix} + 2\left(\alpha + \dfrac{1}{\varepsilon}\right) \begin{pmatrix} I_m & 0 \\ 0 & I_m \end{pmatrix}.$

Now, using the definition of a viscosity subsolution (for $V$) of (5.1) via the superjet (resp. supersolution for $W$ via the subjets), we obtain that, for all $(t, s, x, y) \in [0, T]^2 \times \mathbb{R}^{2m}$ and $(a, p, S) \in \mathcal{P}^{1,2,+} V(t, x)$ [resp. $(b, q, S') \in \mathcal{P}^{1,2,-} W(t, x)$],

$a + \gamma V(t, x)$

$\qquad - \inf_{(\pi, u) \in \overline{U}} \Bigg( pb(x, \pi, u)$

$\qquad\qquad + \dfrac{V(t, x + u\sigma(x, \pi, u)) - V(t, x) - u p \sigma(x, \pi, u)}{u^2} \mathbf{1}_{\{u \neq 0\}}$

$\qquad\qquad + \mathbf{1}_{\{u = 0\}} \langle S\sigma(x, \pi, u), \sigma(x, \pi, u) \rangle \Bigg) \leq 0$

and

$$b + \gamma W(s,y)$$
$$- \inf_{(\pi,u) \in \overline{U}} \Bigg( qb(y,\pi,u)$$
$$+ \mathbf{1}_{\{u \neq 0\}} \frac{W(s, y + u\sigma(y,\pi,u)) - W(s,y) - uq\sigma(y,\pi,u)}{u^2}$$
$$+ \mathbf{1}_{\{u=0\}} \langle S'\sigma(y,\pi,u), \sigma(y,\pi,u) \rangle \Bigg) \geq 0.$$

Indeed, for all $(a,p,S) \in \mathcal{P}^{1,2,+}V(t,x)$, there is some $\varphi \in C^{1,2}([0,T] \times \mathbb{R}^m)$ such that $V - \varphi \leq V(t,x) - \varphi(t,x)$ and $(a,p,S) = (\frac{\partial}{\partial t}\varphi(t,x), \nabla \varphi(t,x), D^2\varphi(t,x))$ (see [12]). Then the above relation follows from Definition 5.1 with $\delta$ less or equal to $\delta_0$ from Remark 6.1. A symmetrical argument gives the above relation for $(b,q,S') \in \mathcal{P}^{1,2,-}W(t,x)$.

Combining the above, we get

$$0 \geq a - b + \gamma(V(t,x) - W(s,y))$$
(6.12)
$$+ \inf_{(\pi,u) \in \overline{U}} \{qb(y,\pi,u) - pb(x,\pi,u)$$
(6.13)
$$+ \mathbf{1}_{\{u \neq 0\}} u^{-2}(W(s, y + u\sigma(y,\pi,u)) - W(s,y)$$
$$- V(t, x + u\sigma(x,\pi,u)) + V(t,x)$$
$$+ u(p\sigma(x,\pi,u) - q\sigma(y,\pi,u)))$$
$$+ \mathbf{1}_{\{u=0\}}(\langle S'\sigma(y,\pi,u), \sigma(y,\pi,u) \rangle$$
$$- \langle S\sigma(x,\pi,u), \sigma(x,\pi,u) \rangle)\}.$$

Now, in light of (6.10), we can find two sequences:

$$\mu_n \triangleq (t_n, x_n, a_n, p_n, S_n)$$
$$\to \left( \hat{t}, \hat{x}, \frac{\hat{t} - \hat{s}}{\varepsilon} + \frac{\alpha}{2} \frac{1}{(T-\hat{t})^2}, \frac{\hat{x} - \hat{y}}{\varepsilon} + \alpha \hat{x}, X \right),$$

$$\nu_n \triangleq (s_n, y_n, b_n, q_n, S'_n)$$
$$\to \left( \hat{s}, \hat{y}, \frac{\hat{t} - \hat{s}}{\varepsilon} - \frac{\alpha}{2} \frac{1}{(T-\hat{s})^2}, \frac{\hat{x} - \hat{y}}{\varepsilon} - \alpha \hat{y}, Y \right),$$

as $n \to \infty$. Here, for all $n \in \mathbb{N}$, $\mu_n = (a_n, p_n, S_n) \in \mathcal{P}^{1,2,+}V(t_n, x_n)$ [resp. $\nu_n = (b_n, q_n, T_n) \in \mathcal{P}^{1,2,-}W(s_n, y_n)$]. We now apply (6.12) to each pair $(\mu_n, \nu_n)$ and then try to take the limit. Note that the special form of $U$, thanks to

(H3), guarantees that at least along a subsequence, again denoted by itself, it holds that

$$\begin{aligned}
0 \geq &\frac{\alpha}{2}\left(\frac{1}{(T-\hat{t})^2} + \frac{1}{(T-\hat{s})^2}\right) + \gamma(V(\hat{t},\hat{x}) - W(\hat{s},\hat{y})) \\
&+ \inf_{(\pi,u)\in\overline{U}}\bigg\{\frac{\hat{x}-\hat{y}}{\varepsilon}(b(\hat{x},\pi,u) - b(\hat{y},\pi,u)) \\
&\quad + \alpha(\hat{x}b(\hat{x},\pi,u) + \hat{y}b(\hat{y},\pi,u)) \\
&\quad + \mathbf{1}_{\{u\neq 0\}}u^{-2}\bigg(W(\hat{s},\hat{y}+u\sigma(\hat{y},\pi,u)) - W(\hat{s},\hat{y}) \\
&\qquad\qquad - V(\hat{t},\hat{x}+u\sigma(\hat{x},\pi,u)) + V(\hat{t},\hat{x}) \\
&\qquad\qquad + u\bigg(\frac{\hat{x}-\hat{y}}{\varepsilon}(\sigma(\hat{x},\pi,u) - \sigma(\hat{y},\pi,u)) \\
&\qquad\qquad\qquad + \alpha(\hat{x}\sigma(\hat{x},\pi,u) + \hat{y}\sigma(\hat{y},\pi,u))\bigg)\bigg) \\
&\quad + \mathbf{1}_{\{u=0\}}(\langle Y\sigma(\hat{y},\pi,u),\sigma(\hat{y},\pi,u)\rangle \\
&\qquad\qquad - \langle X\sigma(\hat{x},\pi,u),\sigma(\hat{x},\pi,u)\rangle)\bigg\}.
\end{aligned}$$

(6.14)

(6.15)

Now, applying (6.11), we have

$$\begin{aligned}
&\langle Y\sigma(\hat{y},\pi,u),\sigma(\hat{y},\pi,u)\rangle - \langle X\sigma(\hat{x},\pi,u),\sigma(\hat{x},\pi,u)\rangle \\
&= \sigma^T(\hat{y},\pi,u)Y\sigma(\hat{y},\pi,u) - \sigma^T(\hat{x},\pi,u)X\sigma(\hat{x},\pi,u) \\
&\geq -\frac{3}{\varepsilon}|\sigma(\hat{x},\pi,u) - \sigma(\hat{y},\pi,u)|^2 \\
&\quad - 2\left(\alpha + \frac{1}{\varepsilon}\right)(|\sigma(\hat{x},\pi,u)|^2 + |\sigma(\hat{y},\pi,u)|^2).
\end{aligned}$$

(6.16)

Recalling the definition of $(\hat{t},\hat{x},\hat{s},\hat{y})$, we have

$$\Psi_{\varepsilon,\alpha}(\hat{t},\hat{x},\hat{s},\hat{y}) \geq \Psi_{\varepsilon,\alpha}(\hat{t},\hat{x}+u\sigma(\hat{x},\pi,u),\hat{s},\hat{y}+u\sigma(\hat{y},\pi,u)).$$

Using relations (6.14) and (6.16), we get, after some simple straightforward calculations,

$$0 \geq \gamma(V(\hat{t},\hat{x}) - W(\hat{s},\hat{y})) - C\alpha(1 + |\hat{x}| + |\hat{y}|) - C\frac{1}{\epsilon}|\hat{x}-\hat{y}|^2,$$

where $C$ is a constant depending only on the bound and Lipschitz constants of $b$ and $\sigma$. Sending $\varepsilon \to 0$ and then $\alpha \to 0$ in the above (taking the "lim sup" if necessary), and noting (6.9), we obtain that

$$\eta \geq \gamma\theta.$$

Since $\eta$ can be arbitrarily small, we see that this is a contradiction. The proof is now complete. $\square$

## APPENDIX: COUNTEREXAMPLE

We now construct an example to show that in general equation (2.7) may not have a unique solution. Let us consider the case $d=1$ and assume that the Lévy measure is of the form $\nu(dx) \stackrel{\triangle}{=} \frac{1}{x^2}\mathbf{1}_{\{x>0\}}$.

Suppose that on some probability space we are given a Brownian motion $B$ and a Poisson random measure $\mu$ whose Lévy measure is $\nu$. We assume that $B$ and $\mu$ are independent, define the stopping time $S = \inf\{t \geq 0, B_t = 1\}$ and set $u_s \stackrel{\triangle}{=} \mathbf{1}_{[S,+\infty)}(s)$. We shall construct two pairs of coefficients $(\alpha_t, \beta_t)$ and $(\alpha'_t, \beta'_t)$, such that the corresponding normal martingales $X$ and $X'$ associated by (2.6) satisfy both (2.7) but are not identical in law.

To this end, let us set $\alpha_t = \mathbf{1}_{[0,S]}(t)$, $\alpha'_t = -\mathbf{1}_{[0,S]}(t)$ and $\beta_t(x) = \beta'_t(x) = \mathbf{1}_{A_t}(x)$, where

$$A_t = \begin{cases} \varnothing, & \text{on } [0, S], \\ [1, \infty), & \text{on } [S, +\infty). \end{cases}$$

Clearly, the process $N_t \stackrel{\triangle}{=} \int_0^t \int_1^{+\infty} \mu(ds\,de) = \mu([0,t] \times [1,\infty))$, $t \geq 0$, is a standard Poisson process independent of $B$. Now, define $\tilde{N}_t \stackrel{\triangle}{=} N_t - t$, $t \geq 0$. Applying Proposition 2.1, we see that both processes

$$X_t = B_{S \wedge t} + \tilde{N}_t - \tilde{N}_{S \wedge t},$$
$$X'_t = -B_{S \wedge t} + \tilde{N}_t - \tilde{N}_{S \wedge t}, \qquad t \geq 0,$$

satisfy the structure equation for the above defined process $u$. We now argue that $X$ and $X'$ are not equal in law. Indeed, we write the stopping time $S$ as $S = \inf\{t \geq 0, X_t = 1\}$ and set $S' = \inf\{t \geq 0, X'_t = 1\}$. Then, it is readily seen that $X \stackrel{\mathcal{D}}{=} X'$ if and only if $(X, S) \stackrel{\mathcal{D}}{=} (X', S')$.

Now consider the stopped processes $X^S_t \stackrel{\triangle}{=} X_{S \wedge t}$ and $(X')^{S'}_t \stackrel{\triangle}{=} (X')_{S' \wedge t}$, $t \geq 0$. If $X \stackrel{\mathcal{D}}{=} X'$, then $X^S \stackrel{\mathcal{D}}{=} (X')^{S'}$ as well. But clearly $X^S$ is continuous, while $(X')^{S'}$ has jumps on the set $\{S < S'\}$, which obviously has a strictly positive probability, a contradiction.

**Acknowledgments.** The authors would like to thank the anonymous referees for their careful readings of the manuscripts and very helpful suggestions, which make many places of the paper more accurate and easier to read. We are indebted to one of the referees for a suggestion that improves the proof of the counterexample significantly.

R. BUCKDAHN
C. RAINER
DÉPARTEMENT DE MATHÉMATIQUES
UNIVERSITÉ DE BRETAGNE OCCIDENTALE
F-29285 BREST CEDEX
FRANCE
E-MAIL: Rainer.Buckdahn@univ-brest.fr
       Catherine.Rainer@univ-brest.fr

J. MA
DEPARTMENT OF MATHEMATICS
UNIVERSITY OF SOUTHERN CALIFORNIA
LOS ANGELES, CALIFORNIA 90089
USA
E-MAIL: jinma@usc.edu